\magnification=1200

\hsize=11.25cm    
\vsize=18cm     
\parindent=12pt   \parskip=5pt     

\hoffset=.5cm   
\voffset=.8cm   

\pretolerance=500 \tolerance=1000  \brokenpenalty=5000

\catcode`\@=11

\font\eightrm=cmr8         \font\eighti=cmmi8
\font\eightsy=cmsy8        \font\eightbf=cmbx8
\font\eighttt=cmtt8        \font\eightit=cmti8
\font\eightsl=cmsl8        \font\sixrm=cmr6
\font\sixi=cmmi6           \font\sixsy=cmsy6
\font\sixbf=cmbx6

\font\tengoth=eufm10 
\font\eightgoth=eufm8  
\font\sevengoth=eufm7      
\font\sixgoth=eufm6        \font\fivegoth=eufm5

\skewchar\eighti='177 \skewchar\sixi='177
\skewchar\eightsy='60 \skewchar\sixsy='60

\newfam\gothfam           \newfam\bboardfam

\def\tenpoint{
  \textfont0=\tenrm \scriptfont0=\sevenrm \scriptscriptfont0=\fiverm
  \def\rm{\fam\z@\tenrm}
  \textfont1=\teni  \scriptfont1=\seveni  \scriptscriptfont1=\fivei
  \def\oldstyle{\fam\@ne\teni}\let\old=\oldstyle
  \textfont2=\tensy \scriptfont2=\sevensy \scriptscriptfont2=\fivesy
  \textfont\gothfam=\tengoth \scriptfont\gothfam=\sevengoth
  \scriptscriptfont\gothfam=\fivegoth
  \def\goth{\fam\gothfam\tengoth}
  
  \textfont\itfam=\tenit
  \def\it{\fam\itfam\tenit}
  \textfont\slfam=\tensl
  \def\sl{\fam\slfam\tensl}
  \textfont\bffam=\tenbf \scriptfont\bffam=\sevenbf
  \scriptscriptfont\bffam=\fivebf
  \def\bf{\fam\bffam\tenbf}
  \textfont\ttfam=\tentt
  \def\tt{\fam\ttfam\tentt}
  \abovedisplayskip=12pt plus 3pt minus 9pt
  \belowdisplayskip=\abovedisplayskip
  \abovedisplayshortskip=0pt plus 3pt
  \belowdisplayshortskip=4pt plus 3pt 
  \smallskipamount=3pt plus 1pt minus 1pt
  \medskipamount=6pt plus 2pt minus 2pt
  \bigskipamount=12pt plus 4pt minus 4pt
  \normalbaselineskip=12pt
  \setbox\strutbox=\hbox{\vrule height8.5pt depth3.5pt width0pt}
  \let\bigf@nt=\tenrm       \let\smallf@nt=\sevenrm
  \normalbaselines\rm}

\def\eightpoint{
  \textfont0=\eightrm \scriptfont0=\sixrm \scriptscriptfont0=\fiverm
  \def\rm{\fam\z@\eightrm}
  \textfont1=\eighti  \scriptfont1=\sixi  \scriptscriptfont1=\fivei
  \def\oldstyle{\fam\@ne\eighti}\let\old=\oldstyle
  \textfont2=\eightsy \scriptfont2=\sixsy \scriptscriptfont2=\fivesy
  \textfont\gothfam=\eightgoth \scriptfont\gothfam=\sixgoth
  \scriptscriptfont\gothfam=\fivegoth
  \def\goth{\fam\gothfam\eightgoth}
  
  \textfont\itfam=\eightit
  \def\it{\fam\itfam\eightit}
  \textfont\slfam=\eightsl
  \def\sl{\fam\slfam\eightsl}
  \textfont\bffam=\eightbf \scriptfont\bffam=\sixbf
  \scriptscriptfont\bffam=\fivebf
  \def\bf{\fam\bffam\eightbf}
  \textfont\ttfam=\eighttt
  \def\tt{\fam\ttfam\eighttt}
  \abovedisplayskip=9pt plus 3pt minus 9pt
  \belowdisplayskip=\abovedisplayskip
  \abovedisplayshortskip=0pt plus 3pt
  \belowdisplayshortskip=3pt plus 3pt 
  \smallskipamount=2pt plus 1pt minus 1pt
  \medskipamount=4pt plus 2pt minus 1pt
  \bigskipamount=9pt plus 3pt minus 3pt
  \normalbaselineskip=9pt
  \setbox\strutbox=\hbox{\vrule height7pt depth2pt width0pt}
  \let\bigf@nt=\eightrm     \let\smallf@nt=\sixrm
  \normalbaselines\rm}

\tenpoint

\def\pc#1{\bigf@nt#1\smallf@nt}         \def\pd#1 {{\pc#1} }

\frenchspacing

\def\raggedbottom{\topskip 10pt plus 36pt\r@ggedbottomtrue}

\def\pointir{\unskip . --- \ignorespaces}

\def\Medbreak{\vskip-\lastskip\medbreak}

\long\def\th#1 #2\enonce#3\endth{
   \Medbreak\noindent
   {\pc#1} {#2\unskip}\pointir{\it #3}\smallskip}

\def\decale#1{\smallbreak\hskip 28pt\llap{#1}\kern 5pt}
\def\decaledecale#1{\smallbreak\hskip 34pt\llap{#1}\kern 5pt}
\def\puce{\smallbreak\hskip 6pt{$\scriptstyle\bullet$}\kern 5pt}

\def\eqalign#1{\null\,\vcenter{\openup\jot\m@th\ialign{
\strut\hfil$\displaystyle{##}$&$\displaystyle{{}##}$\hfil
&&\quad\strut\hfil$\displaystyle{##}$&$\displaystyle{{}##}$\hfil
\crcr#1\crcr}}\,}

\catcode`\@=12

\showboxbreadth=-1  \showboxdepth=-1

\newcount\numerodesection \numerodesection=1
\def\section#1{\bigbreak
 {\bf\number\numerodesection.\ \ #1}\nobreak\medskip
 \advance\numerodesection by1}

\mathcode`A="7041 \mathcode`B="7042 \mathcode`C="7043 \mathcode`D="7044
\mathcode`E="7045 \mathcode`F="7046 \mathcode`G="7047 \mathcode`H="7048
\mathcode`I="7049 \mathcode`J="704A \mathcode`K="704B \mathcode`L="704C
\mathcode`M="704D \mathcode`N="704E \mathcode`O="704F \mathcode`P="7050
\mathcode`Q="7051 \mathcode`R="7052 \mathcode`S="7053 \mathcode`T="7054
\mathcode`U="7055 \mathcode`V="7056 \mathcode`W="7057 \mathcode`X="7058
\mathcode`Y="7059 \mathcode`Z="705A


\def\diagram#1{\def\normalbaselines{\baselineskip=0pt\lineskip=5pt}
\matrix{#1}}

\def\vfl#1#2#3{\llap{$\textstyle #1$}
\left\downarrow\vbox to#3{}\right.\rlap{$\textstyle #2$}}

\def\hfl#1#2#3{\smash{\mathop{\hbox to#3{\rightarrowfill}}\limits
^{\textstyle#1}_{\textstyle#2}}}

\def\Q{{\bf Q}}

\def\R{{\bf R}}
\def\C{{\bf C}}
\def\N{{\bf N}}

\def\U{{\bf U}}
\def\Z{{\bf Z}}

\def\F{{\bf F}}

\def\Hom{\mathop{\rm Hom}\nolimits}

\def\Inf{\mathop{\rm Inf}\nolimits}
\def\Sup{\mathop{\rm Sup}\nolimits}
\def\Card{\mathop{\rm Card}\nolimits}
\def\Gal{\mathop{\rm Gal}\nolimits}
\def\Ker{\mathop{\rm Ker}\nolimits}

\def\to{\rightarrow}

\def\droite#1{\,\hfl{#1}{}{8mm}\,}

\def\mod{\mathop{\rm mod.}\nolimits}
\def\pmod#1{\;(\mod#1)}

\def\sup{\mathop{\rm Sup}\displaylimits}
\def\inf{\mathop{\rm Inf}\displaylimits}

\def\boxit#1{\vbox{\hrule\hbox{\vrule\kern1pt
       \vbox{\kern1pt#1\kern1pt}\kern1pt\vrule}\hrule}}
\def\cqfd{\hfill\boxit{\phantom{\i}}}

\newcount\numero \numero=1
\def\numeroter{{({\oldstyle\number\numero})}\ \advance\numero by1}

\newcount\refno 
\long\def\ref#1:#2<#3>{                                        
\global\advance\refno by1\par\noindent                              
\llap{[{\bf\number\refno}]\ }{#1} \pointir{\it #2} #3\goodbreak }

\newbox\bibbox
\setbox\bibbox\vbox{
\bigskip
\centerline{---$*$---$*$---}
\bigbreak
\centerline{{\pc REFERENCE}}

\ref{\pc TATE} (John):
Problem 9 : The general reciprocity law,
<in Mathematical developments arising from Hilbert problems
(Proc. Sympos. Pure Math., Northern Illinois Univ., De Kalb, Ill.,
1974), pp. 311--322. Proc. Sympos. Pure Math., Vol. XXVIII,
Amer. Math. Soc., Providence, R. I., 1976. > 
\newcount\tate \global\tate=\refno
}



\centerline{\it Eight lectures on qudratic reciprocity}
\bigskip
\centerline{Chandan Singh Dalawat}

\vfill

$$
\vbox{
Si $p$ est numerus primus formae $4n+1$, erit $+p$, si vero $p$
formae $4n + 3$, erit $-p$ residuum vel non residuum cuiusuis numeri
primi qui positive acceptus ipsius $p$ est residuum vel non
residuum.
}
$$
\vskip.2\baselineskip
\rightline{--- Carl Friedrich Gau\ss, {\it Disquisitiones
    arithmeticae}, 1801, \S131.}

\bigskip

$$
\vbox{
Legendre a devin{\'e} la formule et Gauss est devenu instatan{\'e}ment
c{\'e}l{\`e}bre en la prouvant.  En trouver des
g{\'e}n{\'e}ralisations, par exemple aux anneaux d'entiers
alg{\'e}briques, ou d'autres d{\'e}monstrations a constitu{\'e} un
sport national pour la dynastie allemande suscit{\'e} par Gauss
jusqu'{\`a} ce que le reste du monde, {\`a} commencer par le Japonais
Takagi en 1920 et {\`a} continuer par Chevalley une dizaine
d'ann{\'e}es plus tard, d{\'e}couvre le sujet et, apr{\`e}s 1945, le
fasse exploser.  Gouvern{\'e} par un Haut Commissariat qui surveille
rigoureusement l'alignement de ses Grandes Pyramides, c'est
aujourd'hui l'un des domaines les plus respect{\'e}s des
Math{\'e}matiques.  
}
$$
\vskip.2\baselineskip
\rightline{--- Roger Godement, {\it Analyse math{\'e}matique} IV, 2003, p.~313.}

\vfill\eject
\hbox{}
\vfill\eject

\centerline{\bf Lecture 1}
\medskip
\centerline{$\lambda_p(q)=\lambda_q(\lambda_4(p)p)$}
\bigskip

\numeroter The group $\Z^\times$ of invertible elements of the ring
$\Z$ consists of $1$ and $-1$, and will sometimes be thought of as a
($1$-dimensional) vector space over the field $\F_2$ (with the unique
basis $-1$).  

\numeroter Let $G$ be a group.  A {\it character\/} (of $G$) of order
dividing~$2$ is a homomorphism $\chi:G\to\Z^\times$, so an element of
$\Hom(G,\Z^\times)$.  Such a $\chi$ will be called a {\it quadratic
character\/} (of $G$) if it is surjective.

\numeroter The groups of interest will initially be the groups
$G_m=(\Z/m\Z)^\times$ of invertible elements in the rings $\Z/m\Z$ for
$m>0$.  The groups $G_1$ and $G_2$ are trivial.  For $m>2$, we will
often identify $\Z^\times$ with its image in $G_m$.  Recall that if
$\gcd(m,m')=1$, then the canonical map $G_{mm'}\to G_m\times G_{m'}$
is an isomorphism (as a consequence of the ``Chinese remainder
theorem'').

\numeroter Let $p$ be a prime number.  We denote by $\Z_{(p)}$ the
smallest subring of $\Q$ containing $l^{-1}$ for every prime $l\neq
p$.  Let $A$ be a ring.  There is a homomorphism of rings
$f:\Z_{(p)}\to A$ if and only if $l\in A^\times$ for every prime
$l\neq p$.  If $f$ exists, it is unique.

\numeroter Every $a\in\Q^\times$ can be uniquely written as
$a=p^r\alpha$ with $r\in\Z$ and $\alpha\in\Z_{(p)}^\times$.  We define
$v_p(a)=r$, and note that $v_p:\Q^\times\to\Z$ is a surjective
homomorphism satisfying
$
v_p(a+b)\ge\Inf(v_p(a),v_p(b))
$ 
(where $v_p(0)=+\infty$ by convention), with equality if
$v_p(a)\neq v_p(b)$.   We have
$a\in\Z_{(p)}\Leftrightarrow v_p(a)\in\N$ and
$a\in\Z_{(p)}^\times\Leftrightarrow v_p(a)=0$.

\numeroter For every $n>0$, there is a canonical map
$\Z/p^n\Z\to\Z_{(p)}/p^n\Z_{(p)}$ of rings, and the universal property
of $\Z_{(p)}$ provides a map in the other direction, so the two rings
are canonically isomorphic.  The resulting morphism of groups
$\Z_{(p)}^\times\to G_{p^n}$ is surjective for every $n>0$.  Thus,
every $\F_p$-space can be viewed as a $\Z_{(p)}$-module~; in
particular, $\Z^\times$ can be viewed as a $\Z_{(2)}$-module, giving a
meaning to the expression $(-1)^a$ for every $a\in\Z_{(2)}$.

\numeroter Recall that for every prime $p$, the group
$G_p=\F_p^\times$ is {\it cyclic\/} of order $p-1$.  If $p=2p'+1$ is
odd, then $p-1=2p'$ is even, so there is a unique quadratic character
$\lambda_p:\F_p^\times\to\Z^\times$, and $\Ker(\lambda_p)=\F_p^{\times
  2}$.  One can view $\lambda_p$ as a quadratic character of $G_{p^n}$
(for every $n>0$) and of $\Z_{(p)}^\times$ via the surjections of
these groups onto $\F_p^\times$.  We shall see later that $\lambda_p$
is the only quadratic character of $G_{p^n}$ ($n>0$, $p\neq2$).

\numeroter Concretely, for every odd prime $p$ and for every
$a\in\Z_{(p)}^\times$, we have $\lambda_p(a)=+1$ if and only if
$a\equiv x^2\pmod p$ for some $x\in\Z_{(p)}^\times$~; otherwise,
$\lambda_p(a)=-1$.

\numeroter The map $\Z^\times\to(\Z/4\Z)^\times$ is an isomorphism.
We denote the reciprocal isomorphism by
$\lambda_4:(\Z/4\Z)^\times\to\Z^\times$.  We can view $\lambda_4$ as a
quadratic character of $G_{2^n}$ (for every $n>1$) and of
$\Z_{(2)}^\times$ via the surjections of these groups onto
$(\Z/4\Z)^\times$.  Concretely, for every $a\in\Z_{(2)}^\times$, we have 
$$
\lambda_4(a)=(-1)^{a-1\over 2}
=\cases{+1&if $a\equiv+1\pmod4$,\cr
        -1&if $a\equiv-1\pmod4$.\cr}
$$
Note that if $a\in\Z_{(2)}^\times$, then $a-1\in2\Z_{(2)}$, so
$(-1)^{a-1\over 2}$ has a meaning (???).

\numeroter When we view $\Z^\times$ as a subgroup of
$G_8=(\Z/8\Z)^\times$, the quotient $G_8/\Z^\times$ has order~$2$, and
hence it has a unique quadratic character
$\lambda_8:G_8/\Z^\times\to\Z^\times$.  We can view $\lambda_8$ as a
quadratic character of $G_{2^n}$ (for every $n>2$) and of
$\Z_{(2)}^\times$ via the surjections of these groups onto $G_8$ and
thence onto $G_8/\Z^\times$.  It can be easily checked that, for every
$a\in\Z_{(2)}^\times$,
$$
\lambda_8(a)=(-1)^{a^2-1\over 8}
=\cases{+1&if $a\equiv\pm1\pmod8$,\cr
        -1&if $a\equiv\pm5\pmod8$.\cr}
$$ 
Note that if $a\in\Z_{(2)}^\times$, then $a^2-1\in8\Z_{(2)}$, so
$(-1)^{a^2-1\over 8}$ has a meaning (???).  When $a=2a'+1$ for some
$a'\in\Z_{(2)}$, the definition amounts to
$$
\lambda_8(a)
=\cases{(-1)^{a'}&if $a'\equiv0,1\pmod4$,\cr
        (-1)^{a'-1}&if $a'\equiv2,3\pmod4$.\cr}
$$
It can also be easily verified that
$\lambda_8(a)=(-1)^{\lambda_4(a)a-1\over4}$ for every $a\in\Z_{(2)}^\times$.

\numeroter Notice that the only quadratic characters of
$G_8=(\Z/8\Z)^\times$ are $\lambda_4$, $\lambda_8$ and
$\lambda_4\lambda_8$.  Of these, only $\lambda_4$ comes from a
quadratic characters of $G_4=(\Z/4\Z)^\times$, and only $\lambda_8$ is
even in the sense that $\lambda_8(-1)=+1$.  Note that
$\lambda_4,\lambda_8$ is a basis of the $\F_2$-space
$\Hom(G_8,\Z^\times)$.  We shall see later that for every $n>2$, the
only quadratic characters of $G_{2^n}$ are $\lambda_4$, $\lambda_8$
and $\lambda_4\lambda_8$.  

\numeroter Incidentally, the unique quadratic character
$\lambda_\infty$ of $\R^\times$ is given by
$\lambda_\infty(a)=a/|a|_\infty$, where $|a|_\infty=\Sup(a,-a)$, and
notice that $a\in\R^\times$ is a square if and only if
$\lambda_\infty(a)=+1$, if and only if $a>0$.  In view of
$\R^\times=\Z^\times\times\R^\times_+$, the uniqueness of
$\lambda_\infty$ follows from the fact that every character
$\chi:\R^\times_+\to\Z^\times$ is trivial~: $\chi(a)=\chi(\sqrt
a)^2=1$.  The character $\lambda_\infty$ is sometimes denoted
$\mathop{\rm sgn}\nolimits$ (for the Latin {\it signum\/}).

\numeroter  The group $\Z_{(p)}^\times$ is generated by $-1$ and the
set of primes $l$ distinct from $p$, so any homomorphism
$\lambda:\Z_{(p)}^\times\to H$ ($H$ being a group) is uniquely
determined by $\lambda(-1)$ and the $\lambda(l)$.  When $p=2$, we have
explicit formul\ae\ for $\lambda_4$ and $\lambda_8$ (and hence also
for $\lambda_4\lambda_8$, namely
$$
\lambda_4\lambda_8(a)=(-1)^{{a-1\over4}+{a^2-1\over8}}
$$
for every $a\in\Z_{(2)}^\times$).  When $p$ is an odd prime,
$\lambda_p$ is completely determined by $\lambda_p(-1)$,
$\lambda_p(2)$ and the $\lambda_p(q)$ for every odd primes $q\neq p$.

\numeroter Let $p=2p'+1$ be an odd prime.  The {\it quadratic
  reciprocity law\/} asserts that 
$$
\lambda_p(-1)
=\lambda_4(p),\quad 
\lambda_p(2)
=\lambda_8(p),\quad\hbox{and\ }
\lambda_p(q)
=\lambda_q(\lambda_4(p)p)
$$ 
for every odd prime $q=2q'+1$ distinct from $p$.
It was discovered by Euler and independently by Legendre, who gave a
partial proof, and it was finally proved by Gau\ss\ at the age of 18.
Gau\ss\ called it the {\it theorema fundamentale\/} and gave at least
five, six, seven, or eight different proofs --- the count depending on
the historian consulted --- during the course of his life.  We shall
give below one of the simplest proofs.

\numeroter The quotation from the {\it Disquitiones arithmeticae\/}
(\S131) of Gau\ss\
$$
\vbox{\it Si $p$ est numerus primus formae $4n+1$, erit $+p$, si vero $p$
formae $4n + 3$, erit $-p$ residuum vel non residuum cuiusuis numeri
primi qui positive acceptus ipsius $p$ est residuum vel non
residuum.}
$$
can be translated into our notation as
$\lambda_p(q)=\lambda_q(\lambda_4(p)p)$. 

\numeroter (Euler) {\it Let\/ $p=2p'+1$ be an odd prime.  We have\/
  $\lambda_p(a)=a^{p'}$ for every\/ $a\in\F_p^\times$.  Equivalently,
  $\lambda_p(a)\equiv a^{p'}\pmod p$ for every\/
  $a\in\Z_{(p)}^\times$.}

{\it Proof}.  Let $a\in\F_p^\times$, and let $r\in\F_p^\times$ be a
generator (an element of order $p-1=2p'$, or a ``\thinspace primitive
root modulo~$p$\thinspace'')~; we have $r^{p'}=-1$ (because
$r^{2p'}=1$ and $r^{p'}\neq1$).  Write $a=r^n$ for some $n\in\Z$
(uniquely determined $\mod{2p'}$, and hence $\mod2$) and note that
$\lambda_p(a)=+1$ if and only if $n\equiv0\pmod2$.  Also,
$a^{p'}=r^{np'}=(-1)^n$, which equals $+1$ if and only if
$n\equiv0\pmod2$.  Hence $\lambda_p(a)=a^{p'}$. \cqfd

\numeroter Applying Euler's criterion (??) to $a=-1$ we get the
supplementary law $\lambda_p(-1)=(-1)^{p'}=\lambda_4(p)$.  Another
easy consequence is that for every odd prime $q=2q'+1$,
$$
\lambda_q(\lambda_4(p))
=\lambda_q((-1)^{p'})
=(-1)^{p'q'},
$$
which explains why the law $\lambda_p(q)=\lambda_q(\lambda_4(p)p)$
(for $q$ distinct from $p$) is often written as
$\lambda_p(q)=(-1)^{p'q'}\lambda_q(p)$.

\numeroter  A character $\chi:G_m\to\Z^\times$ will be called
{\it even\/} if $\chi(-1)=1$, {\it odd\/} if $\chi(-1)=-1$.  Among the
quadrtic character we have named, $\lambda_8$ is even, $\lambda_4$ and
$\lambda_4\lambda_8$ are odd, and, for an odd prime $p$, the character
$\lambda_p$ is even if $\lambda_4(p)=1$, odd if $\lambda_4(p)=-1$
(??).  We also agree to call $\lambda_\infty$ (??) odd because
$\lambda_\infty(-1)=-1$.

\numeroter We say that a subset $S\subset\F_p^\times$ is a {\it
  section\/} if the canonical projection modulo $\Z^\times$ induces a
bijection $S\to\F_p^\times/\Z^\times$ of sets.  Every section has $p'$
elements, and the map $(e,x)\mapsto ex$ is a bijection
$\Z^\times\times S\to\F_p^\times$.  Every $a\in\F_p^\times$ uniquely
determines a map $e_a:S\to\Z^\times$ and a permutation $\sigma_a:S\to
S$ such that $ax=e_a(x)\sigma_a(x)$ for every $x\in S$.  The simplest
section is $S=[1,p']$, and there are $2^{p'}$ sections in all.

\numeroter (Gau\ss) {\it Let\/ $S\subset\F_p^\times$ be a section
  (??).  We have\/ $\lambda_p(a)=\prod_{x\in S}e_a(x)$ for every\/
  $a\in\F_p^\times$.}

{\it Proof}.  It is sufficient (??) to prove that $a^{p'}=\prod_{x\in
  S}e_a(x)$.  Multiplying together the relations
$ax=e_a(x)\sigma_a(x)$ (for $x\in S$), we get
$$
a^{p'}\prod_{x\in S}x
=\prod_{x\in S}e_a(x)\sigma_a(x)
=\prod_{x\in S}e_a(x).\prod_{x\in S}\sigma_a(x)
=\prod_{x\in S}e_a(x).\prod_{x\in S}x
$$
and hence the result, because $\prod_{x\in S}x$ is invertible in
$\F_p$.  \cqfd

\numeroter Let us prove the second supplementary law
$\lambda_p(2)=\lambda_8(p)$.  Take the section $S=[1,p']$ and take
$a=2$ in (??).  If $p'=4n$ or $p'=4n+1$ for some $n>0$, then, for all
$x\in S$,
$$
e_2(x)=\cases{+1&if $x\in[1,2n]$,\cr
 -1&if $x\notin[1,2n]$,\cr}
$$
so $\prod_{x\in S}e_2(x)=(-1)^{p'-2n}=(-1)^{p'}$.
Similarly, if $p'=4n+2$ or $p'=4n+3$ for
some $n>0$, then
$$
e_2(x)=\cases{+1&if $x\in[1,2n+1]$,\cr
 -1&if $x\notin[1,2n+1]$,\cr}
$$
so $\prod_{x\in S}e_2(x)=(-1)^{p'-(2n+1)}=(-1)^{p'-1}$.  We have seen
that 
$$
\lambda_p(2)
=\cases{(-1)^{p'}&if $p'\equiv0,1\pmod4$,\cr
        (-1)^{p'-1}&if $p'\equiv2,3\pmod4$.\cr}
$$ 
Comparison with  (??) leads to the desired result~:
$\lambda_p(2)=\lambda_8(p)$. \cqfd

\numeroter Let us prove that $\lambda_p(q)=(-1)^{p'q'}\lambda_q(p)$
for any two distinct odd primes $p=2p'+1$ and $q=2q'+1$.  The idea is
to compute the product of all elements in
$(\F_p^\times\times\F_q^\times)/\Z^\times$ (where
$\Z^\times\subset(\F_p^\times\times\F_q^\times)$ is embedded
diagonally) in two different ways, by taking two different systems of
representatives.  This proof is inspired by (??) and was given by
Rousseau in 1991.  Another proof, directly based upon (??), is given
in the exercises~; it is due to Frobenius in 1914, and goes back to
the fifth proof of Gau\ss.

\numeroter One system of representatives is
$[1,p-1]\times[1,q']\subset(\F_p^\times\times\F_q^\times)$.  So the
representatives are
$$
\matrix{
(1,1),&(1,2),&\cdots,&(1,q'),\cr
(2,1),&(2,2),&\cdots,&(2,q'),\cr
\cdots,&\cdots,&\cdots,&\cdots,\cr
(p-1,1),&(p-1,2),&\cdots,&(p-1,q'),\cr
}
$$ 
and their product is visibly $((p-1)!^{q'},q'!^{p-1})$.  Notices
that $(p-1)!=-1$ in $\F_p^\times$ (``Wilson's theorem'', which can be
proved by taking the product of all elements in $\F_p^\times$), and
that
$$
q'!^{p-1}
=(q'!^2)^{p'}
=((-1)^{q'}(q-1)!)^{p'}
=((-1)^{q'}(-1))^{p'}
=(-1)^{p'q'+p'}
$$
in $\F_q^\times$, so the product of all the representatives is
$((-1)^{q'},(-1)^{p'q'+p'})$, which is equivalent (under $\Z^\times$)
to 
$$
(1,(-1)^{p'q'+p'+q'}).
$$

\numeroter Secondly, note that for every integer $m>0$, the set of
integers in $[1,m]$ which are prime to $n=2m+1$ is a system of
representatives in $(\Z/n\Z)^\times$ for
$(\Z/n\Z)^\times\!/\Z^\times$.  Take $m=q'p+p'=p'q+q'$, so that $n=pq$
and we have an isomorphism of groups
$(\Z/n\Z)^\times\to\F_p^\times\times\F_q^\times$ (??) inducing the
identity on the subgroups $\Z^\times$.  This gives our second system
of representatives~: the set of integers in $[1,m]$ which are prime to
$pq$.

\numeroter Their product in $\F_p^\times$ is computed by first
computing the product of all integers in $[1,m]$ which are prime to
$p$, namely
$$
\matrix{
1,&2,&\cdots,&(p-1),\cr
p+1,&p+2,&\cdots,&p+(p-1),\cr
\cdots,&\cdots,&\cdots,&\cdots,\cr
(q'-1)p+1,&(q'-1)p+2,&\cdots,&(q'-1)p+(p-1),\cr
}
$$
and $q'p+1,q'p+2,\cdots,q'p+p'$, and then dividing by the product of
all multiples of $q$ in $[1,m]$.  The product, in $\F_p^\times$, of
all these prime-to-$p$ elements of $[1,m]$ is
$(p-1)!^{q'}p'!=(-1)^{q'}p'!$.  Also, the multiples of $q$ in $[1,m]$
are
$$
1.q,\ 2.q,\ \cdots,\ p'.q
$$
and their product in $\F_p^\times$ is $p'!q^{p'}=p'!\lambda_p(q)$, so
the product, in $\F_p^\times$, of all integers in $[1,m]$ which are
prime to $pq$ is
$$
{(-1)^{q'}p'!\over p'!\lambda_p(q)}
=(-1)^{q'}\lambda_p(q).
$$
Similarly, the product, in $\F_q^\times$, of all integers in $[1,m]$
which are prime to $pq$ is $(-1)^{p'}\lambda_q(p)$.  

\numeroter So the product, in $\F_p^\times\times\F_q^\times$, of this
system of representatives (of
$(\F_p^\times\times\F_q^\times)\!/\Z^\times$ in
$\F_p^\times\times\F_q^\times$, namely the set of integers in $[1,m]$
which are prime to $pq=2m+1$) is $((-1)^{q'}\lambda_p(q),
(-1)^{p'}\lambda_q(p))$, which is equivalent (under $\Z^\times$) to
$$ 
(1,(-1)^{p'+q'}\lambda_p(q)\lambda_q(p)).
$$ 
Comparing this result with (??) gives
$\lambda_p(q)\lambda_q(p)=(-1)^{p'q'}$.  \cqfd

\bigbreak
\leftline{\it Exercises}
\medskip

\numeroter Let $p=2p'+1$ be an odd prime.  In the lemma (??), take $S$
to be section $[1,p']\subset\F_p^\times$ and take $a=q$, where
$q=2q'+1$ is an odd prime distinct from $p$, to conclude that
$\lambda_p(q)=(-1)^M$, where $M$ is the number of pairs
$(x,y)\in[1,p']\times[1,q']$ such that $qx-py\in[-p',-1]$.  ({\it
  Hint~:} The relation $e_q(x)=-1$ is equivalent to the existence of a
$y\in\Z$ such that $qx-py\in[-p',-1]$.  This $y$, if it exists, is
unique.  Show that $y\in[1,q']$.)


\numeroter Similarly prove that $\lambda_q(p)=(-1)^N$, where $N$ is
the number of pairs $(x,y)\in[1,p']\times[1,q']$ such that
$qx-py\in[1,q']$ and hence $\lambda_p(q)\lambda_q(p)=(-1)^{M+N}$.
Show that there is no pair $(x,y)\in[1,p']\times[1,q']$ such that
$qx-py=0$, so the exponent $M+N$ is also the number of pairs such that
$qx-py\in[-p',q']$.


\numeroter Show that the map $(x,y)\mapsto(p'+1-x,q'+1-y)$ is an
involution of $[1,p']\times[1,q']$ inducing a bijection between the
set of pairs $(x,y)$ satisfying $qx-py<-p'$ with the set of pairs
$(x,y)$ satisfying $qx-py>q'$, using the identity
$$
q(p'+1-x)-p(q'+1-y)=-(qx-py)-p'+q'.
$$ 

\numeroter Conclude that $M+N\equiv p'q'\pmod2$, thereby proving the
law $\lambda_p(q)\lambda_q(p)=(-1)^{p'q'}$.



\medbreak

\numeroter (Bost, 2012).  It is known that
$p=2^{43112609}-1$ is a prime number~; it was in fact the largest
prime known in the year $2012$.  Compute $\lambda_p(2012)$.  ({\it
  Hint :\/} If a prime $p$ and an $a\in\F_p^\times$ are fixed, then
$a^n$ depends only on $n\pmod{p-1}$, and can be computed by repeated
squarings if we know the base-$2$ expansion of $n$.)

\numeroter We have the prime decomposition $2012=2^2.503$, therefore
$\lambda_p(2012)=\lambda_p(503)$.  Show that
$2^{43112609}\equiv92\pmod{503}$, and hence $p\equiv91\pmod{503}$.
({\it Hint :\/} As $43112609\equiv347\pmod{502}$, we have
$2^{43112609}\equiv2^{347}\pmod{503}$, and as
$347=1+2+2^3+2^4+2^6+2^8$, we just need to compute $2^{2^i}\pmod{503}$
by repeated squarings for $i\in[1,8]$ and multiply some of them
together to get $2^{347}\equiv92\pmod{503}$. I thank Paul Vojta for
correcting a mistake in an earlier version of this calculation.)


\numeroter Show that $\lambda_p(503) =\lambda_{503}(-p)
=\lambda_{503}(-91)$.  Now use the prime decomposition $91=7.13$ to
complete exercise (??).

\vfill\eject

\centerline{\bf Lecture 2}
\medskip
\centerline{$\prod_v(a,b)_v=1$}
\bigskip

\numeroter For each prime number $p$, we shall define a bilinear map
(the quadratic hilbertian symbol at $p$)
$$
(\ ,\ )_p:\Q^\times\times\Q^\times\to\Z^\times
$$ 
using the quadratic character $\lambda_p$ (??) for $p\neq2$ and the
quadratic characters $\lambda_4$ (??), $\lambda_8$ (??) for $p=2$.  We
shall also define the (quadratic) hilbertian symbol $(\ ,\ )_\infty$
at the ``archimedean place'' $\infty$ (using the quadratic character
$\lambda_\infty$), and show that the three laws (??)
$$\lambda_p(-1)=\lambda_4(p),\quad 
\lambda_p(2)=\lambda_8(p),\quad
\lambda_p(q)=\lambda_q(\lambda_4(p)p)
$$
(in which $p$ and $q$ are distinct odd primes) can be encapsulated in
a single neat product formula $\prod_v(a,b)_v=1$, where $v$ runs over
all {\it places\/} of $\Q$, namely the prime numbers and also the
archimedean place $\infty$, and $a,b\in\Q^\times$.

\numeroter First a bit of notation.  The quadratic characters
$\lambda_*$ (for $*$ an odd prime or $*=4,8,\infty$) take values in
the multiplicative group $\Z^\times$.  Sometimes we need to think of
them as taking values in the field $\F_2$ (over which $\Z^\times$ is a
$1$-dimensional vector space, with basis $-1$), and then we denote
them by $\varepsilon_*$, so the relation between the two is
$$
\lambda_*(a)=(-1)^{\varepsilon_*(a)}
$$ 
for every $a\in\Z_{(p)}^\times$ (where $p=2$ for $*=4,8$ and $p=*$ if
$*$ is an odd prime) or for every $a\in\R^\times$ if $*=\infty$.  

\numeroter The advantage of this notation is that we can replace
conditions such as ``\thinspace $\lambda_4(a)=1$ or
$\lambda_4(b)=1$\thinspace'' by the condition ``\thinspace
$\varepsilon_4(a)\varepsilon_4(b)=0$\thinspace''.  Similarly, ``$a>0$
or $b>0$\thinspace'' is equivalent to ``\thinspace
$\varepsilon_\infty(a)\varepsilon_\infty(b)=0$\thinspace''.  The
property $\lambda_*(ab)=\lambda_*(a)\lambda_*(b)$ implies that
$\varepsilon_*(ab)=\varepsilon_*(a)+\varepsilon_*(b)$.

\numeroter Let us make the definition explicit.  For every
$a\in\R^\times$, we have
$$
\varepsilon_\infty(a)=\cases{0&if $a>0$,\cr
                            1&if $a<0$.\cr}
$$

\numeroter Write $x\equiv_ny$ for $x\equiv y\pmod n$.  For
$a\in\Z_{(2)}^\times$, we have (??)
$$
\varepsilon_4(a)\equiv_2{a-1\over2}
 \equiv_2 \cases{0&if $a\equiv_4+1$,\cr
                 1&if $a\equiv_4-1$;\cr}
$$ 
note the curiostiy
$\varepsilon_\infty(a)=\varepsilon_4(\lambda_\infty(a))$ for every
$a\in\R^\times$, which has its uses.  Similarly, for
$a\in\Z_{(2)}^\times$, 
$$
\varepsilon_8(a)\equiv_2{a^2-1\over8}
 \equiv_2 \cases{0&if $a\equiv_8\pm1$,\cr
                 1&if $a\equiv_8\pm5$.\cr}
$$
by (??).  As $\lambda_8$ is a morphism of groups (??), we have
${ab-1\over2}\equiv_2{a-1\over2}+{b-1\over2}$ and
${(ab)^2-1\over8}\equiv_2{a^2-1\over8}+{b^2-1\over8}$ for
$a,b\in\Z_{(2)}^\times$.

\numeroter Finally, for every odd prime $p$ and every
$a\in\Z_{(p)}^\times$, we have
$$
\varepsilon_p(a)=
\cases{0&if $\bar a\in\F_p^{\times2}$\cr
       1&if $\bar a\notin\F_p^{\times 2}$.\cr}
$$

\numeroter The quadratic reciprocity law (???) can of course be rewritten as
$$
\varepsilon_p(-1)=\varepsilon_4(p),\quad \varepsilon_p(2)=\varepsilon_8(p),\quad
\varepsilon_p(q)=\varepsilon_4(p)\varepsilon_4(q)+\varepsilon_q(p)
$$ 
(for any two distinct odd primes $p$ and $q$) but this reformulation
is no improvement. Reformulating the law in terms of hilbertian
symbols is going to be a substantial improvement.

\numeroter Put
$(a,b)_\infty=(-1)^{\varepsilon_\infty(a)\varepsilon_\infty(b)}$, so
that $(a,b)_\infty=1$ if and only if $a>0$ or $b>0$, if and only if
there exist $x,y\in\R$ such that $ax^2+by^2=1$.

\numeroter Now let $p$ be a prime number.  Note that every
$x\in\Q^\times$ can be uniquely written as $x=p^{v_p(x)}u_x$, with
$v_p(x)\in\Z$ and $u_x\in\Z_{(p)}^\times$ (???).  Let
$a,b\in\Q^\times$, write
$$
a=p^{v_p(a)}u_a,\quad
b=p^{v_p(b)}u_b,\quad 
(v_p(a),v_p(b)\in\Z, u_a,u_b\in\Z_{(p)}^\times)~;
$$ 
and note that $u_{ab}=u_au_b$ in addition to $v_p(ab)=v_p(a)+v_p(b)$.  Put
$$
t_{a,b}
=(-1)^{v_p(a)v_p(b)}a^{v_p(b)}b^{-v_p(a)}
=(-1)^{v_p(a)v_p(b)}u_a^{v_p(b)}u_b^{-v_p(a)}
$$ 
which is visibly in $\Z_{(p)}^\times$ (??).  Notice that
$t_{b,a}=t_{a,b}^{-1}$ and $t_{a,bc}=t_{a,b}t_{a,c}$ (for every
$c\in\Q^\times$).

\numeroter If $p\neq2$, define
$
(a,b)_p
=(-1)^{\varepsilon_p(t_{a,b})} 
$, 
so that $(a,b)_p=+1$ if and only if
$\bar t_{a,b}\in\F_p^{\times2}$~;
otherwise, $(a,b)_p=-1$.
It follows from the definitions that
$$
\eqalign{
(a,b)_p
=\lambda_p(t_{a,b})
&=(-1)^{\varepsilon_4(p)v_p(a)v_p(b)}\lambda_p(u_a)^{v_p(b)}\lambda_p(u_b)^{-v_p(a)}\cr
&=(-1)^{\varepsilon_4(p)v_p(a)v_p(b)+v_p(b)\varepsilon_p(u_a)-v_p(a)\varepsilon_p(u_b)}.\cr
}
$$

\numeroter For $p=2$, define
$
(a,b)_2
=(-1)^{\varepsilon_4(u_a)\varepsilon_4(u_b)+\varepsilon_8(t_{a,b})}
$,
so that $(a,b)_2=+1$ if and only if
${u_a-1\over2}{u_b-1\over2}+{t_{a,b}^2-1\over8}\equiv0\pmod2$~;
otherwise, $(a,b)_2=-1$.  Again, it follows from the definitions that
$$
(a,b)_2
=(-1)^{{u_a-1\over2}{u_b-1\over2}+{t_{a,b}^2-1\over8}}
=(-1)^{{u_a-1\over2}{u_b-1\over2}}\lambda_8(t_{a,b}).
$$
Notice finally that $\varepsilon_8(t_{a,b})\equiv
v_2(b)\varepsilon_8(u_a)-v_2(a)\varepsilon_8(u_b)\pmod2$, therefore
$$
(a,b)_2=(-1)^{\varepsilon_4(u_a)\varepsilon_4(u_b)+v_2(b)\varepsilon_8(u_a)-v_2(a)\varepsilon_8(u_b)}.
$$

\numeroter  These definitions might seem a bit contrived.  Once we have
introduced the fields $\Q_p$ (containing $\Q$) for primes $p$, we will
see that $(a,b)_p=1$ if and only if there exist $x,y\in\Q_p$ such that
$ax^2+by^2=1$.  This is the natural definition of $(a,b)_p$, valid for
all $a,b\in\Q_p^\times$, and the above formul\ae\ are the result of a
computation.  The natural definition brings out the analogy with the
symbol $(\ ,\ )_\infty$ (??) on $\R=\Q_\infty$.

\numeroter The hilbertian symbol $(\ ,\ )_v$ (where $v$ is a prime $p$
or $v=\infty$) possesses some elementary properties which we enumerate
next.  Most of them are straightforward calculations, and completely
obvious if $v=\infty$.  So assume that $v=p$ is a prime, and
$a=p^{v_p(a)}u_a,b=p^{v_p(b)}u_b$ as above.

\numeroter $(a,b)_v=(b,a)_v$.  (Interchanging $a,b$ replaces
$t_{a,b}$ by $t_{a,b}^{-1}$ (and interchanges $v_p(a),
v_p(b)$ and $u_a,u_b$), but
$\lambda_*(t_{a,b}^{-1})=\lambda_*(t_{a,b})$ for $*=p$ if
$p\neq2$ and $*=8$ if $p=2$.)

\numeroter $(a,bc)_v=(a,b)_v(a,c)_v$.  (Indeed,
$v_p(bc)=v_p(b)+v_p(c)$, $u_{bc}=u_bu_c$ and
$t_{a,bc}=t_{a,b}t_{b,c}$, and the $\varepsilon_*$ are
homomorphisms (for $*=p$ if $p\neq2$ and $*=4,8$ if $p=2$).)

\numeroter It follows from these two facts that $(a,b)_v$ depends only
on the classes of $a,b$ modulo $\Q^{\times 2}$ in the sense that
$(a,b)_v=(ac^2,b)_v$ for every $c\in\Q^\times$ and
$(a,b)_v=(a,bd^2)_v$ for every $d\in\Q^\times$ (so we may assume
that $a,b$ are squarefree integers).

\numeroter $(a,-a)_v=1$.  (For $t_{a,b}=1$ (and
$\varepsilon_4(u_{-a})=1+\varepsilon_4(u_a)$ when $p=2$).)

The following two propositions are immediate consequences of the
definitions.

\numeroter {\it Let\/ $p$ be an odd prime and let\/ $a,b\in\Z_{(p)}^\times$.  Then}
$$
(a,b)_p=1,\quad (a,pb)_p=(-1)^{\varepsilon_p(a)}=\lambda_p(a).
$$ \cqfd

\numeroter {\it Let\/ $a,b\in\Z_{(2)}^\times$.  Then we have\/}
$$
(a,b)_2=(-1)^{\varepsilon_4(a)\varepsilon_4(b)},\quad
(a,2)_2=(-1)^{\varepsilon_8(a)}=\lambda_8(a).
$$ \cqfd


\numeroter As an example, we have $(p,p)_p=(-p,p)_p(-1,p)_p=(-1,p)_p$,
which is $=\lambda_p(-1)=\lambda_4(p)$ if $p\neq2$ (??), and
$=\lambda_8(-1)=1$ if $p=2$ (??).  Note that this example and (??)
(for $p\neq2$) or (??) (for $p=2$) cover all possible cases, by
symmetry (??) and multiplicativity (??).

\numeroter {\it Let\/ $a,b\in\Q^\times$.  Then\/ $(a,b)_v=1$ for
  almost every place\/ $v$ of\/ $\Q$, and we have the product formula}
$$
\prod_v(a,b)_v=1.
$$

{\it Proof}.  That $(a,b)_v=1$ for almost every $v$ follows from (??)
and the fact that $a,b\in\Z_{(p)}^\times$ for almost every prime $p$.
By the symmetry and bilinearity of the symbol $(\ ,\ )_v$ and the fact
that the $\F_2$-space $\Q^\times\!/\Q^{\times 2}$ admits a basis
consisting of $-1$ and the primes numbers, it is sufficient to prove
the product formula in the following cases~:

\numeroter {\it $a=p$ and\/ $b=q$, where\/ $p$ and $q$ are distinct
  odd primes.}  We then have $(p,q)_v=1$ for all $v\neq p,q,2$, and
$$
(p,q)_p=\lambda_p(q),\quad 
(p,q)_q=\lambda_q(p),\quad 
(p,q)_2=(-1)^{\varepsilon_4(p)\varepsilon_4(q)},
$$ 
by (??) and (??), so the product formula follows in this case from the
law $\lambda_p(q)\lambda_q(p)=(-1)^{\varepsilon_4(p)\varepsilon_4(q)}$
(??).

\numeroter {\it $a=p$, where $p$ is an odd prime, and\/ $b=2$.}
Similarly, we have $(p,2)_v=1$ for all $v\neq p,2$, and
$$
(p,2)_p=\lambda_p(2),\quad 
(p,2)_2=\lambda_8(p)
$$
so the product formula follows from the law
$\lambda_p(2)=\lambda_8(p)$ (??).

\numeroter {\it $a=p$, where $p$ is an odd prime, and\/ $b=-1$.}  We
then have $(p,-1)_v=1$ for all $v\neq p,2$, and
$$
(p,-1)_p=\lambda_p(-1),\quad 
(p,-1)_2=(-1)^{\varepsilon_4(p)\varepsilon_4(-1)}
=(-1)^{\varepsilon_4(p)}
$$
so the product formula follows from the law
$\lambda_p(-1)=\lambda_4(p)$ (??). 

\numeroter {\it $a=2$ and\/ $b=-1$.}  We then have $(2,-1)_v=1$ for
all $v\neq 2,\infty$, and
$$
(2,-1)_2=\lambda_8(-1)=1,\quad
(2,-1)_\infty=1,
$$
so the product formula holds trivially in this case.

\numeroter {\it $a=-1$ and\/ $b=-1$.}  We then have $(-1,-1)_v=1$ for
$v\neq2,\infty$ and $(-1,-1)_v=-1$ for $v=2,\infty$.  

We don't need to consider the case $(a,b)=(l,l)$ for any prime $l$
because $(l,l)_v=(-1,l)_v(-l,l)_v=(-1,l)_v$ by (??), which has been
treated in (??) if $l\neq2$ and in (??) if $l=2$.  This completes the
proof of the product formula (??) in all cases.  \cqfd

\numeroter  Conversely, it is obvious that the product formula (???)
implies the quadratic reciprocity law (???).

\numeroter The above proof is summarised in the following table (in
which $l$ is a prime $\neq2,p,q$ and the blank entries stand for $+1$)
$$
\vbox{\openup1\jot
\halign{\hfil$#$\qquad\qquad&$#$\qquad&$#$\qquad&$#$\qquad&$#$\qquad&$#$\cr
v&\infty&2&p&q&l\cr
\noalign{\vskip-3\jot}
\multispan{5}\hrulefill\cr
(p,q)_v&&(-1)^{\varepsilon_4(p)\varepsilon_4(q)}&\lambda_p(q)&\lambda_q(p)\cr
(p,2)_v&&\lambda_8(p)&\lambda_p(2)\cr
(p,-1)_v&&\lambda_4(p)&\lambda_p(-1)\cr
(2,-1)_v\cr 
(-1,-1)_v&-1&-1\cr
\noalign{\vskip-3\jot}
\multispan{5}\hrulefill.\cr
}}
$$

\numeroter {\it Fix\/ $a,b\in\Q^\times$.  The number of places $v$
  such that\/ $(a,b)_v=-1$ is (finite and) even.}   \cqfd

\numeroter We shall see later that given any finite set $S$ of places
of $\Q$ such that $\Card S$ is even, there exist $a,b\in\Q^\times$
such that $(a,b)_v=-1$ for $v\in S$ and $(a,b)_v=1$ for $v\notin S$.

\numeroter {\it Let\/ $w$ be a place of\/ $\Q$.  If\/ $(a,b)_v=1$ for
  all\/ $v\neq w$, then\/ $(a,b)_w=1$.}  \cqfd

\numeroter For example, if $(a,b)_p=1$ for every prime number $p$,
then at least one of $a,b$ must be $>0$.

\bigbreak
\leftline{\it Exercises}
\medskip

\numeroter Let $a=\lambda_\infty(a)\prod_{p\neq2}p^{v_p(a)}$ (??) be
the prime decomposition of an {\it odd\/} integer $a\in\Z$, so that
$\lambda_\infty(a)$ is the sign of $a$ and $v_p(a)=0$ for almost all
primes~$p$, and define
$\psi_a(n)
=\prod_{v_p(a)\equiv1\pmod2}\lambda_p(n)$
for every integer $n$ prime to $a$.  Show that
$$
\psi_a(-1)=(-1)^{\varepsilon_4(a)+\varepsilon_\infty(a)},\qquad
\psi_a(2)=(-1)^{\varepsilon_8(a)},
$$
and, for every odd integer $b\in\Z$ prime to $a$, we have the
reciprocity law
$$
\psi_a(b)
=(-1)^{\varepsilon_4(a)\varepsilon_4(b)+\varepsilon_\infty(a)\varepsilon_\infty(b)}
 \psi_b(a).
$$

\numeroter Define $k(a)=\prod_{v_p(a)\equiv1\pmod2}p$, and view
$\psi_a$ as a character of $G_{k(a)}$ (??) via the isomorphism
$G_{k(a)}\to\prod_{v_p(a)\equiv1\pmod2}\F_p^\times$ (??). 
Show that the above reciprocity law continues to remain valid for any
two odd integers $a,b\in\Z$ such that $\gcd(a,k(b))=1$ and
$\gcd(b,k(a))=1$ (or in other words $\bar a\in G_{k(b)}$ and $\bar
b\in G_{k(a)}$).

\numeroter Let $a\in\Z$ be a squarefree integer, and put
$m=4|a|_\infty$.  Show that there is a unique homomorphism
$\chi_a:G_m\to\C^\times$ such that $\chi_a(p)=\lambda_p(a)$ for every
prime~$p$ not dividing~$m$.  Moreover, $\chi_a^2=1$, but $\chi_a\neq1$
if $a\neq1$.  ({\it Hint~:\/} It is clear that $\chi_a$ is unique (if
it exists) and has order dividing~$2$.  As for the existence, take
$\chi_b=\lambda_4^{\varepsilon_4(b)}\lambda_{l_1}\lambda_{l_2}\ldots\lambda_{l_r}$
if $b=l_1l_2\ldots l_r$ (where $r=0$ if $b=1$) is a product of
distinct odd primes $l_i$, and take
$$
\chi_{-b}=\lambda_4\chi_b,\quad
\chi_{2b}=\lambda_8\chi_b,\quad
\chi_{-2b}=\lambda_4\lambda_8\chi_b~;
$$ 
this defines $\chi_a$ for every (squarefree) $a$.  Suppose that
$a\neq1$.  If $a=-1,2$ or $-2$, then clearly $\chi_a\neq1$.  Otherwise
$a$ has some odd prime factors $l_1,l_2,\ldots,l_r$ ($r>0$)~; for any
$x\in\Z$ such that $\lambda_{l_1}(x)=-1$ and $x\equiv1\pmod{4l_2\ldots
  l_r}$, we have $\chi_a(x)=-1$.)

\numeroter Show that $\chi_a(x)=\prod_{l\mid
  m}(a,x)_l=\prod_{\gcd(l,m)=1}(a,x)_l$ for every integer $x>0$ prime
to~$m$.


\vfill\eject

\centerline{\bf Lecture 3}
\medskip
\centerline{$\Z_p$}
\bigskip

\numeroter Let $p$ be a prime number.  For every $n>0$, we have the
finite ring $A_n=\Z/p^n\Z$ with $p^n$ elements, and a surjective
morphism of rings $\varphi_n:A_{n+1}\to A_n$, with kernel
$p^nA_{n+1}$, so that we have the exact sequence $0\to A_1\to
A_{n+1}\to A_n\to0$, where the first map is ``\thinspace
multiplication by $p^n$\thinspace'' (the unique morphism of groups
$A_1\to A_{n+1}$ such that $1\mapsto p^n$).  Similarly we have the
exact sequence $0\to A_n\to A_{n+1}\to A_1\to0$, where the first map
is ``\thinspace multiplication by $p$\thinspace''.

\numeroter A $p$-{\it adic integer\/} is a system of elements
$(x_n)_{n>0}$ such that $x_n\in A_n$ and $\varphi_n(x_{n+1})=x_n$.
The set of $p$-adic integers is denoted by $\Z_p$~; it is a subset of
the product $\prod_{n>0}A_n$ defined by the vanishing of
$\varphi_m\circ\pi_{m+1}-\pi_m$ for all $m>0$, where $\pi_m$ is the
natural projection $\prod_{n>0}A_n\to A_m$.  The restriction of
$\pi_m$ to $\Z_p$ is surjective because the $\varphi_n$ are
surjective.

\numeroter For every $n>0$, base-$p$ expansion in $\N$ gives a natural
bijection $[0,p[^n\to[0,p^n[$, namely
        $(b_i)_{i\in[0,n[}\mapsto\sum_{i\in[0,n[}b_ip^i$, and thence a
                natural bijection $[0,p[^n\to A_n$.  If $x_{n+1}\in
                    A_{n+1}$ corresponds to $(b_i)_{i\in[0,n]}$, then
                    $\varphi_n(x_{n+1})$ corresponds to
                    $(b_i)_{i\in[0,n[}$.  It follows that the set
                        $\Z_p$ is in natural bijection with the
                        product $[0,p[^{\N}$~; a $p$-adic integer
                            $x\in\Z_p$ corresponds to a sequence
                            $(b_i)_{i\in\N}$ if and only if
$$
\pi_n(x)\equiv\sum_{i\in[0,n[}b_ip^i\pmod{p^n}
$$
(in $A_n$) for every $n>0$.  In particular, the set $\Z_p$ has the
cardinality of the continuum.

\numeroter Let $x=(x_n)_{n>0}$ be a $p$-adic integer.  If $x_n=0$
(resp.~$x_n=1$) for every $n>0$, then we write $x=0$ (resp.~$x=1$).
If $y=(y_n)_{n>0}$ is another $p$-adic integer, we define
$$
-x=(-x_n)_{n>0},\quad
x+y=(x_n+y_n)_{n>0},\quad
xy=(x_ny_n)_{n>0}.
$$ 
These definitions give $\Z_p$ the structure of a commutative ring for
which each $\pi_m:\Z_p\to A_m$ is a morphism of rings.

\numeroter For every $m>0$ define $V_m=\pi_m^{-1}(0)$.  There is a
unique topology on $\Z_p$ for which $(x+V_m)_{m>0}$ is a fundamental
system of open neighbourhoods of $x$, for every $x\in\Z_p$.  Each
$x+V_m$ is also closed in $\Z_p$ because $A_m$ is finite.  This
topology is compatible with the ring structure of $\Z_p$, and each
$\pi_m$ is continuous.  The space $\Z_p$ is compact because it is a
closed subset of the product $\prod_{n>0}A_n$.

\numeroter In other words, {\it the profinite ring\/ $\Z_p$ is the
  projective limit of the inverse system\/} $(\varphi_n:A_{n+1}\to
A_n)_{n>0}$.  This means that given any ring $A$ and homomorphisms of
rings $f_n:A\to A_n$ such that $\varphi_n\circ f_{n+1}=f_n$ for every
$n>0$, there is a unique homomorphism of rings $\iota:A\to\Z_p$ such
that $f_n=\pi_n\circ\iota$ for every $n>0$.  Indeed,
$\iota(a)=(f_n(a))_{n>0}$ for every $a\in A$.  Also, $X$ being a
space, a map $f:X\to\Z_p$ is continuous if and only if $\pi_n\circ f$
is continuous for every $n>0$.

\numeroter If we consider instead the inverse system
$(\varphi_n:B_{n+1}\to B_n)_{n>0}$, where $B_n=\F_p[T]/(T^n)$ and $T$
is an indeterminate, we get the profinite ring $\F_p[[T]]$ which is
similar to $\Z_p$ in many respects.

\numeroter Let $\iota:\Z\to\Z_p$ be the natural morphism of rings.  For every
$m>0$, the composite $\pi_m\circ\iota$ is the canonical projection
$\Z\to A_m$ (passage to the quotient modulo $p^m\Z$).  In particular,
$\iota$ is injective, for if $a\in\Z$ is such that $a\equiv
0\pmod{p^m}$ for every $m>0$, then $a=0$.  We indentify $\Z$ with its
image $\iota(\Z)$ in $\Z_p$.

\numeroter {\it For every\/ $m>0$, multiplication by $p^m$ is
  injective on $\Z_p$, and the ideal\/ $p^m\Z_p$ is the kernel\/ $V_m$
  of\/ $\pi_m:\Z_p\to A_m$.}

{\it Proof}.  For the first part, it suffices to prove that the map
$x\mapsto px$ is injective on $\Z_p$.  Indeed, if $px=0$ for some
$x=(x_n)_{n>0}$ in $\Z_p$, then $px_{n+1}=0$ for every $n>0$, and
there exist $y_{n+1}\in A_{n+1}$ such that $x_{n+1}=p^ny_{n+1}$.  But
then
$$
x_n=\varphi_n(x_{n+1})=p^n\varphi_n(y_{n+1})=0
$$ 
for every $n>0$, and hence $x=0$. It follows that multiplication by
$p^m$ is injective on $\Z_p$ for every $m>0$. 

Clearly $p^m\Z_p\subset V_m$.  If $x\in V_m$, then $x_m=0$, so
$x_{m+r}\in p^mA_{m+r}$ for every $r>0$, and $x_{m+r}=p^my_r$ for some
$y_r\in A_{m+r}$ uniquely determined $\pmod{p^r}$ and such that
$\varphi_r(y_{r+1})=y_r$.  The $p$-adic integer $y=(y_r)_{r>0}$ is
such that $x=p^my$, and therefore $p^m\Z_p=V_m=\Ker\pi_m$. \cqfd

\numeroter We thus have the exact sequence $0\to
p^m\Z_p\to\Z_p\to\Z/p^m\Z\to0$ for every $m>0$~; in particular,
$\F_p=\Z_p/p\Z_p$.

\numeroter {\it For every\/ $m>0$, the induced map\/
  $\Z/p^m\Z\to\Z_p/p^m\Z_p$ is an isomorphism of rings.} \cqfd

\numeroter {\it For\/ $x\in\Z_p$, the following conditions are
  equivalent~:
$$
(1)\quad x\in\Z_p^\times,\qquad (2)\quad x\notin p\Z_p,\qquad 
(3)\quad\pi_1(x)\in\F_p^\times.
$$
In particular, $p\Z_p$ is the unique maximal ideal of the ring\/
$\Z_p$.}

{\it Proof}.  The equivalence of $(2)$ and $(3)$ follows from the fact
that $\F_p=\Z_p/p\Z_p$.  As the implication $(1)\Longrightarrow(3)$ is
clear, it suffices to prove that $(2)\Longrightarrow(1)$.

Suppose that $x\notin p\Z_p$, so that $\pi_1(x)\in\F_p^\times$.  It
follows that for every $n>0$, we have $x_n\notin pA_n$.  Therefore
there exist $y_n,z_n\in A_n$ such that $x_ny_n=1-pz_n$, or
equivalently $x_nx_n'=1$, with
$$
x_n'=y_n(1+pz_n+\ldots+p^{n-1}z_n^{n-1}).
$$ 
We have $\varphi_n(x_{n+1}')=x_n'$ (because $\varphi_n(x_{n+1}')$ is
also an inverse of $x_n$ in the ring $A_n$), so we get a $p$-adic
integer $x'\in\Z_p$ such that $xx'=1$, and hence
$x\in\Z_p^\times$. \cqfd

\numeroter The group $\Z_p^\times$ can be thought of as the projective
limit of the inverse system $(A_{n+1}^\times\to A_n^\times)_{n>0}$.

\numeroter {\it There is a unique morphism of rings\/
  $\iota:\Z_{(p)}\to\Z_p$.  It is injective, and the composite
  $\pi_m\circ\iota$ is the natural projection $\pmod{p^m\Z_{(p)}}$.}

{\it Proof}.  Indeed, if $u\in S$, where $S\subset\Z$ is the
multiplicative subset of integers prime to $p$, then $u\in\Z_p^\times$
(??), so by the universal property of the localisation
$\Z_{(p)}=S^{-1}\Z$, there exists a unique morphism of rings $\iota:
\Z_{(p)}\to\Z_p$ extending the inclusion $\Z\subset\Z_p$.  It is
injective because the only $x\in\Z_{(p)}$ such that
$x\equiv0\pmod{p^n}$ for every $n>0$ is $x=0$.  Finally, the induced
map $\Z_{(p)}/p^m\Z_{(p)}\to\Z_p/p^m\Z_p$ is an isomorphism because
the composite $\Z/p^m\Z\to\Z_{(p)}/p^m\Z_{(p)}\to\Z_p/p^m\Z_p$ is an
isomorphism (??).  \cqfd

\numeroter {\it Every\/ $x\neq0$ in\/ $\Z_p$ can be uniquely written
  as\/ $x=p^mu$, with $m\in\N$, $u\in\Z_p^\times$.}

{\it Proof}.  Let $x\in\Z_p$.  If $x\neq0$, there is a largest integer
$m>0$ such that $\pi_m(x)=0$~; we then have $x=p^mu$ with $u\notin
p\Z_p$ (??).  By (??), $u\in\Z_p^\times$.  The decomposition $x=p^mu$ is
unique because $m$ is uniquely determined by $x$, and because
$y\mapsto p^my$ is injective (??).  \cqfd

\numeroter For $x=p^mu$ ($m\in\N, u\in\Z_p^\times$), we put
$v_p(x)=m$, and define $v_p(0)=+\infty$.  Note that each of the three
conditions in (??)  is equivalent to
``\thinspace$v_p(x)=0$\thinspace''.  It is clear that
$v_p(xy)=v_p(x)+v_p(y)$.  This definition is compatible with (??), and
the inequality $v_p(x+y)\ge\Inf(v_p(x),v_p(y))$, with equality if
$v_p(x)\neq v_p(y)$, continues to hold for all $x,y\in\Z_p$.

\numeroter {\it The ring\/ $\Z_p$ is integral and every ideal\/
  ${\goth a}\neq0$ is generated by\/ $p^n$ for some\/ $n\in\N$.}

{\it Proof}.  For $x\neq0$ and $y\neq0$ in $\Z_p$, we have
$v_p(x)+v_p(y)<+\infty$ and hence $xy\neq0$. Next, let ${\goth
  a}\neq0$ be an ideal of $\Z_p$, and let $n$ be the smallest number
in $v_p({\goth a})$.  We claim that ${\goth a}=p^n\Z_p$.  First, if
$x\in{\goth o}$ is such that $v_p(x)=n$, we have $x=p^n\alpha$ for
some $\alpha\in\Z_p^\times$, or equivalently $p^n=x.\alpha^{-1}$, so
$p^n\in{\goth a}$.  Secondly, for every $y\neq0$ in $\goth a$, we have
$v_p(y)\ge n$, and $y=p^{v_p(y)-n}\beta.p^n$ for some
$\beta\in\Z_p^\times$ (??), so ${\goth a}\subset p^n\Z_p$.  Hence
${\goth a}=p^n\Z_p$.  \cqfd

\numeroter We put
$|x|_p=p^{-v_p(x)}$ for $x\neq0$ in $\Z_p$, and define $|0|_p=0$.  We
then have $|x-y|_p\le\Sup(|x|_p,|y|_p)$ (with equality if
$|x|_p\neq|y|_p$) (??)  and hence $d_p(x,y)=|x-y|_p$ is a distance on
$\Z_p$ satisfying
$$
d_p(x,z)\le\Sup(d_p(x,y), d_p(y,z)),
$$
({\it the ultrametric inequality\/}), stronger than the
triangular inequality. 

\numeroter {\it The topology on\/ $\Z_p$ can be defined by the
  distance\/ $d_p$, for which it is complete.}

{\it Proof}.  That the topology can be defined by $d_p$ follows from
the fact that the fundamental system of open neighbourhoods $p^n\Z_p$
of $0$ (??) is an open ball for $d_p$, namely $d_p(x,0)< p^{-(n-1)}$.
That $\Z_p$ is complete for $d_p$ is a consequence of its
compactness. \cqfd

\numeroter We can reverse the process and define $\Z_p$ as the
completion of $\Z$ for the distance $d_p$.  Notice that every open
ball in $\Z_p$ is also a closed ball, for example $|x|_p<1$ is the
same as $|x|_p\le p^{-1}$.  Also, any point of a ball can be
considered as its ``centre'', as follows from the ultrametric
inequality.

\numeroter Recall (??) that to every $p$-adic integer $x\in\Z_p$, we
have associated a sequence $(b_i)_{i\in\N}$ of elements $b_i\in[0,p[$
    characterised by the fact that for every $n\in\N$,
$$
\sum_{i\in[0,n]}b_ip^i\equiv\pi_{n+1}(x)\pmod{p^{n+1}}.
$$

\numeroter {\it For every\/ $x\in\Z_p$, the associated series\/
  $\sum_{i\in\N}b_ip^i$ converges in\/ $\Z_p$ to\/ $x$.}

{\it Proof}.  
For every $n\in\N$, let $s_n=\sum_{i\in[0,n]}b_ip^i$ be the partial
sums, and fix an integer $m>0$.  We have to show that almost all $s_n$
are in $x+p^m\Z_p$.  This is clearly the case as soon as $n>m$, for
$s_n-x\in p^{n+1}\Z_p$, by the defining property of the sequence
$(b_i)_{i\in\N}$.  \cqfd

\numeroter As an example, take $x=-1$, so that $b_i=p-1$ for every
$i\in\N$, and hence
$$
-1=\sum_{i\in\N}(p-1)p^i,\quad 
(1-p)^{-1}=\sum_{i\in\N}p^i,
$$
in $\Z_p$.  Both expressions give $-1=1+2+2^2+\ldots$ for $p=2$.

\numeroter {\it The subset\/ $\N$ is dense in\/ $\Z_p$.  More
  generally, if $b\in\Z$ is prime to~$p$ and if $a\in\Z$, then $a+b\N$
  is dense in $\Z_p$}

{\it Proof}.  We have to show that for every $x\in\Z_p$ and every
$n>0$, there exists an $x'\in\N$ such that $x'\in x+V_n$~; it suffices
to take an $x'$ whose image in $A_n$ is $x_n$. The second statement
follows from this because $x\mapsto a+bx$ is an of $\Z_p$ whenever
$\gcd(b,p)=1$, as $|b|_p=1$ isometry (??).  \cqfd

\numeroter So a $p$-adic integer $x$ can be considered as a formal
expression $x=\sum_{i\in\N}b_ip^i$, with $b_i\in[0,p[$.  Addition and
    multiplication can be defined by interpreting the partial sums as
    elements of $\N$ and taking base-$p$ expansions of the sum or
    product.  For $x\neq0$, the valuation $v_p(x)$ is the smallest
    index $i$ such that $b_i\neq0$.  The greater the valuation of $x$,
    the closer $x$ is to $0$ in the $p$-adic sense of $d_p(0,x)$.  For
    example, the sequence $1,p,p^2,\ldots$ converges to $0$ in $\Z_p$.

\bigbreak
\leftline{\it Exercises}
\medskip

\numeroter Let $n>0$ be an integer and let $n=a_\nu
p^\nu+a_{\nu+1}p^{\nu+1}+\cdots$ be its base-$p$ expansion, where
$\nu=v_p(n)$, $a_i\in[0,p[$, $a_\nu>0$, and $a_i=0$ for almost all
    $i$.  Put $s_n=a_\nu+a_{\nu+1}+\cdots$ and
    $t_n=a_\nu!a_{\nu+1}!\ldots$ (with the convention $0!=1$).  Show
    that $p-1$ divides $n-s_n$, that $\bar t_n\in\F_p^\times$, and
    that
$$ v_p(n!)={n-s_n\over p-1},\qquad {n!\over (-p)^{v_p(n!)}}\equiv
    t_n\pmod p.  
$$
Observe that $v_p(n!)=\sum_{j>0}\left\lfloor{n\over
      p^j}\right\rfloor$. ({\it Hint~:} Use induction on $n$,  noting
that $n!=(n-1)!.n$ and that the base-$p$ expansion of $n-1$ is 
$$
n-1=(p-1)+(p-1)p+\cdots+(p-1)p^{\nu-1}+(a_\nu-1)p^\nu+a_{\nu+1}p^{\nu+1}+\cdots.) 
$$ 

\numeroter Let $x\in\Z_2$.  For every integer $n>0$, put
$c_n=\prod_{i\in[0,n[}(1-2i)$.  Show that the series
    $\sum_{n>0}c_n{(4x)^n\over n!}$ converges to some $y\in\Z_2$ and
    that $(1+y)^2=1+8x$.  ({\it Hint~:\/} $2^n$ does not divide $n!$.)
    We shall see later that $1+y$ is the unique square root of $1+8x$
    such that $y\equiv0\pmod4$.

\vfill\eject

\centerline{\bf Lecture 4}
\medskip
\centerline{$\Z_p^\times$}
\bigskip

\numeroter Let us move on to serious things.  Suppose we want to find
the roots in $\Z_p$ of some polynomial $f\in\Z_p[T]$.  This amounts to
finding, for every $n>0$, a root $\xi_n\in A_n$ of $f$ such that
$\xi_{n+1}\equiv\xi_n\pmod{p^n}$, so a first necessary condition for
$f$ to have a root in $\Z_p$ is that there should exist an
$x_1\in\Z_p$ such that $f(x_1)\equiv0\pmod p$.

\numeroter Let us try to improve such an $x_1$ to an $x_2=x_1+pz_2$
(with $z_2\in\Z_p$, so that $x_2\equiv x_1\pmod p$) such that
$f(x_2)\equiv0\pmod{p^2}$.  To compute $f(x_1+pz_2)$ we use the
polynomial identity
$$
f(T+S)=f(T)+f'(T)S+g(T,S)S^2 
$$
(valid for some $g\in\Z_p[T,S]$, where $f'$ denotes the formal
derivative of $f$), so $f(x_2)\equiv f(x_1)+f'(x_1)z_2p\pmod{p^2}$.
As $f(x_1)=y_1p$ for some $y_1\in\Z_p$, 
$$
f(x_1+z_2p)\equiv0\pmod{p^2}\quad\Leftrightarrow\quad
y_1+f'(x_1)z_2\equiv0\pmod p.
$$
So {\it if\/} $f'(x_1)\not\equiv0\pmod p$ (so that
$f'(x_1)\in\Z_p^\times$), we take $z_2=-y_1/f'(x_1)$, and then
$$
x_2=x_1-{f(x_1)\over f'(x_1)},\quad f(x_2)\equiv0\pmod{p^2},\quad
x_2\equiv x_1\pmod p.
$$ 
The moral of this story is that if $f$ has a {\it simple\/} root in
$A_1$, then it can be uniquely lifted to a root of $f$ in $A_2$.  

\numeroter This process can be iterated if it turns out that
$f'(x_2)\in\Z_p^\times$ (which we will see is the case)~: if we put
$x_3=x_2-f(x_2)/f'(x_2)$, then
$$
x_3\equiv x_2\pmod{p^2},\qquad 
f(x_3)\equiv0\pmod{p^3}.
$$
Such is the basic idea behind the proof of the following slightly more
general and important result, known as Hensel's lemma.

\numeroter {\it Let\/ $f\in\Z_p[T]$ and\/ $x\in\Z_p$ be such that
  $f'(x)\neq0$, and put\/ $\delta=v_p(f'(x))$. Suppose that we have\/
  $f(x)\equiv0\pmod{p^m}$ for some $m>2\delta$.
Then there exists a unique\/ $\xi\in\Z_p$ such that\/ 
$$
f(\xi)=0,\quad \xi\equiv x\pmod{p^{m-\delta}},\quad v_p(f'(\xi))=\delta.
$$}

{\it Proof}.  Note first that $f'(x)=p^\delta u$ for some
$u\in\Z_p^\times$ (??), that $f(x)/f'(x)$ is in $\Z_p$, and that
therefore so is $y=x-f(x)/f'(x)$.  We claim that
$$
f(y)\equiv0\pmod{p^{m+1}},\quad y\equiv x\pmod{p^{m-\delta}},\quad
v_p(f'(y))=\delta.
$$
Indeed, write $f(x)=p^ma$ for some $a\in\Z_p$, so that $y-x\in
p^{m-\delta}\Z_p$.  The polynomial identity $(??)$ implies that
$$
f(y)=f(x)-{f(x)\over f'(x)}f'(x)+(y-x)^2t=(y-x)^2t
$$
for some $t\in\Z_p$, which implies that $f(y)\in p^{2m-2\delta}\Z_p$.
But $m>2\delta$ by hypothesis, so $f(y)\equiv 0\pmod{p^{m+1}}$.
Next, applying the polynomial identity $(??)$ to $f'$, we get
$$
f'(y)=f'(x+(y-x))=f'(x)+(y-x)s
$$ 
for some $s\in\Z_p$. Note that $v_p(y-x)>\delta$, so
$v_p(f'(y))=v_p(f'(x))=\delta$ (??), and our claim about $y$ is
established. 

\numeroter Let us come to the proof of (??).  The existence of $\xi$
follows immediately from the preceding discussion.  Indeed, starting
with the given $x_0=x$, this algorithm furnishes an $x_1=y$ to which
the algorithm can be reapplied.  We thus get a sequence $(x_i)_{i\in\N}$
of $p$-adic integers such that
$$
f(x_i)\equiv 0\pmod{p^{m+i}},\quad
x_{i+1}\equiv x_i\pmod{p^{m+i-\delta}},\quad
v_p(f'(x_i))=\delta,
$$
for every $i\in\N$.  The sequence $(x_i)_{i\in\N}$ converges in $\Z_p$
to the desired $\xi$.

\numeroter Let us show the uniqueness $\xi$ (satisfying the stated
conditions).
In fact, we shall show that if $f(\eta)=0$ for some $\eta\in\Z_p$
satisfying the weaker congruence $\eta\equiv x\pmod{p^{\delta+1}}$,
then $\eta=\xi$.  The polynomial identity $(??)$ gives
$$
f(\eta)=f(\xi)+(\eta-\xi).f'(\xi)+(\eta-\xi)^2.a
$$
for some $a\in\Z_p$.  As $\eta$ and $\xi$ are roots of $f$, the above
relation implies that 
$$
(\eta-\xi)(f'(\xi)+(\eta-\xi)a)=0.
$$ 
But because $v_p(f'(\xi))=\delta$ whereas $v_p((\eta-\xi)a)>\delta$
by hypothesis, the second factor is $\neq0$.  The only possibility is
that $\eta=\xi$ and we are done.  \cqfd

\bigbreak

\numeroter We extend the quadratic characters $\lambda_p$ for $p\neq2$
(resp.~$\lambda_4$ and $\lambda_8$ for $p=2$) to the whole of
$\Z_p^\times$ by posing $\lambda_p(u)=\lambda_p(\pi_1(u))$ for
$p\neq2$ (resp.~$\lambda_4(u)=\lambda_4(\pi_2(u))$ and
$\lambda_8(u)=\lambda_8(\pi_3(u))$ for $p=2$), where $\pi_m$ is the
projection $\Z_p^\times\to(\Z/p^m\Z)^\times$.

\numeroter {\it For primes\/ $p\neq2$, a unit $u\in\Z_p^\times$ is a square if and
only if\/ $\lambda_p(u)=+1$.}

{\it Proof}.  Clearly $\lambda_p(x^2)=+1$ for every $x\in\Z_p^\times$.
Suppose that $\lambda_p(u)=+1$.  Then there exists an $x\in\Z_p$ such
that $f(x)\equiv0\pmod p$, where $f=T^2-u$.  As $f'(x)\not\equiv0\pmod
p$, we can apply (??) with $\delta=0$, $m=1$, to conclude that there
is a unique $\xi\in\Z_p$ such that $f(\xi)=0$, $\xi\equiv x\pmod p$,
and $v_p(f'(\xi))=0$ (which just means that $\xi\in\Z_p^\times$).  In
other words, $\xi$ is the square root of $u$ congruent $\pmod p$ to
the given $x$. \cqfd

\numeroter {\it For odd primes $p$, the $\F_2$-space
  $\Z_p^\times/\Z_p^{\times 2}$ consists of $\{\bar 1, \bar u\}$,
  where $u$ is any unit such that $\lambda_p(u)=-1$.}  

{\it Proof}.  That there is such a unit $u\in\Z_p^\times$ follows from
  the fact that the projection $\pi_1:\Z_p^\times\to\F_p^\times$ is
  surjective.  \cqfd

\numeroter {\it For\/ $u\in\Z_2^\times$ to be a square, it is
  necessary and sufficient that\/ $\lambda_4(u)=+1$ and\/
  $\lambda_8(u)=+1$, or equivalently\/ $u\equiv1\pmod8$. }

{\it Proof}.  The equivalence of the conditions ``\thinspace
$\lambda_4(u)=+1$ and\/ $\lambda_8(u)=+1$\thinspace'' and ``\thinspace
$u\equiv1\pmod8$\thinspace'' is easy to see, as it is to see that they
hold for every $u=x^2$ ($x\in\Z_2^\times)$.  If they hold for some
$u\in\Z_2^\times$, then $f=T^2-u$ has the root $x=1\pmod 8$.  We have
$f'(1)=2$, so we can apply (??) with $\delta=1$, $m=3$, to conclude
that there is a unique $\xi\in\Z_2$ such that $\xi^2=u$,
$\xi\equiv1\pmod4$ and $v_2(2\xi)=1$ (which just means that
$\xi\in\Z_2^\times$).  \cqfd

\numeroter {\it The\/ $\F_2$-space\/ $\Z_2^\times/\Z_2^{\times
    2}=(\Z/8\Z)^\times$ has a basis consisting of\/ $\bar5, -\bar1$.
  The values of $\lambda_4,\lambda_8$ on this basis are given by the
  matrix}
$$
\bordermatrix{&\lambda_4&\lambda_8\cr
\phantom{-}\bar5&\phantom{+}1&-1\cr
-\bar1&-1&\phantom{+}1\cr
}.$$ \cqfd

\numeroter Next we prove that there is canonical section
$\omega:\F_p\to\Z_p$ of the projection $\pi_1:\Z_p\to\F_p$.  Of course
$\omega$ cannot be a morphism of groups, much less a morphism of
rings, but it has the desirable property of being {\it
  multiplicative}.

\numeroter {\it Let\/ $R\subset\Z_p$ be the set of roots of\/ $T^p-T$.
  The reduction map\/ $\pi_1:\Z_p\to\F_p$ gives a bijection\/
  $R\to\F_p$.  Moreover, if\/ $x,y\in R$, then\/ $xy\in R$ and if\/
  $\omega:\F_p\to R$ denotes the reciprocal bijection, then
  $\omega(ab)=\omega(a)\omega(b)$.}

{\it Proof}.  Put\/ $f=T^p-T$~; every $a\in\F_p$ is a simple root of
$f\pmod p$, so there is a unique root $\omega(a)\in R$ of $f$ in
$\Z_p$ such that $\pi_1(\omega(a))=a$.  As $R$ (being the set of roots
in an integral ring (??) such as $\Z_p$ of a polynomial of degree $p$)
can have at most $p$ elements, the map $\omega:\F_p\to R$ is
bijective, and $\pi_1$ induces the reciprocal bijection.

If $x,y\in R$, then $x^p=x$ and $y^p=y$, therefore $(xy)^p=xy$, and
hence $xy\in R$.  Finally, the multiplicativity
$\omega(ab)=\omega(a)\omega(b)$ for $a,b\in\F_p$ follows from the
unicity of the $\omega(ab)\in R$ such that $\pi_1(\omega(ab))=ab$ and
the fact that for the element $\omega(a)\omega(b)\in R$ we have
$\pi_1(\omega(a)\omega(b))=ab$. \cqfd 

\numeroter The subset $R$ is called the set of multiplicative
representatives of $\F_p$.  Let $R^\times$ be the set of roots of
$T^{p-1}-1$ in $\Z_p$~; we sometimes identify $\F_p^\times$ with
$R^\times$.  The morphism of groups $\omega:\F_p^\times\to\Z_p^\times$
is a section of the short exact sequence
$$
1\to U_1\to\Z_p^\times\to\F_p^\times\to1.
$$
We shall see later that the torsion subgroup of
$\Q_p^\times$ is $R^\times$ (resp.~$\Z^\times$) for $p\neq2$
(resp.~$p=2$).

\numeroter Let's introduce some notation.  For every $n>0$, let
$U_n=1+p^n\Z_p$ be the group of units in $\Z_p$ which are
$\equiv1\pmod{p^n}$.  The $U_n$ form a decreasing sequence of open
subgroups of $U=\Z_p^\times$, and we have already remarked that $U$
can be identified with the projective limit of the system
$(\varphi_n:U/U_{n+1}\to U/U_n)_{n>0}$.

\numeroter {\it 
For every prime\/ $p$, the group\/ $\Z_p^\times$ is the internal
direct product of\/ $R^\times$ and\/ $U_1$.}

{\it Proof}.  It is clear that $R^\times\cap U_1=\{1\}$, and every
$x\in\Z_p^\times$ can be written as $x=wu$, with $w=\omega(\pi_1(x))$
in $R^\times$ and $u=xw^{-1}$ in $U_1$. \cqfd

\numeroter For every $n>0$, the map $(1+p^nx)\mapsto\pi_1(x)$ defines
upon passage to the quotient an isomorphism of groups
$U_n/U_{n+1}\to\Z/p\Z$, as follows from the identity
$$
(1+p^nx)(1+p^ny)\equiv 1+p^n(x+y)\pmod{p^{n+1}}.
$$
In short, the map $x\mapsto(x-1)/p^n\pmod p$ is a surjective morphism
of groups $U_n\to\F_p$, and its kernel is $U_{n+1}$.

\numeroter {\it For every prime\/ $p$ and every\/ $r\in[1,p[$, the binomial
coefficient\/ ${p \choose r}=p!/r!(p-r)!$ is divisible by\/ $p$.}

{\it Proof}.  Indeed, $v_p(p!)=1$ whereas $v_p(r!)=0$ and
$v_p((p-r)!)=0$.  \cqfd

\numeroter {\it Let\/ $p\neq2$ (resp.~$p=2$) be a prime, and let\/
  $n>0$ (reps.~$n>1$) be an integer.  If\/ $x\in U_n$ but\/ $x\notin
  U_{n+1}$, then\/ $x^p\in U_{n+1}$ but\/ $x^p\notin U_{n+2}$.  In
  other words, $(\ )^p$ induces an isomorphism\/ $U_n/U_{n+1}\to
  U_{n+1}/U_{n+2}$ of groups.}

{\it Proof}.  Write $x=1+p^na$, so that $a\not\equiv0\pmod p$, by
hypothesis.  The binomial theorem gives
$$
x^p=1+p^{n+1}a+\ldots+p^{np}a^p,
$$ 
where the suppressed terms ${p\choose r}(p^na)^r$ ($1<r<p$) are all
divisible by $p^{2n+1}$ (??) and hence also by $p^{n+2}$.  At the same
time, we have $np>n+1$ (because $n>1$ if $p=2$), so we get
$x^p\equiv1+p^{n+1}a\pmod{p^{n+2}}$, which implies that $x^p\in
U_{n+1}$ but $x^p\notin U_{n+2}$.  The induced morphism
$U_n/U_{n+1}\to U_{n+1}/U_{n+2}$ is an isomorphism because it is not
trivial and the two groups are of order~$p$.  \cqfd

\numeroter In short, for every $n>0$ if $p\neq2$ (resp.~$n>1$ if
$p=2$), the map $(\ )^p$ takes $U_n$ to $U_{n+1}$, the composite map
$U_n\to U_{n+1}/U_{n+2}$ is surjective, its kernel in $U_{n+1}$, and
we have a commutative diagram of groups
$$
\diagram{
U_n&\droite{(\ )^p}&U_{n+1}\cr
\vfl{}{}{5mm}&&\vfl{}{}{5mm}\cr
U_n/U_{n+1}&\droite{\sim}&U_{n+1}/U_{n+2}.\cr
} 
$$
  
\numeroter {\it For every\/ $n>0$ if\/ $p\neq2$ (resp.~$n>1$ if\/ $p=2$), the
group\/ $U_1/U_n$ (resp.~$U_2/U_n$) is cyclic of order\/ $p^{n-1}$
(resp.~$2^{n-2})$.}

{\it Proof}.  Let us first treat the case $p\neq2$.  Choose an
$\alpha$ (for example $\alpha=1+p$) such that $\alpha\in U_1$ but
$\alpha\notin U_2$.  By repeated application of (??), we see that
$\alpha^{p^i}\in U_{i+1}$ but $\alpha^{p^i}\notin U_{i+2}$ for every
$i>0$.  In other words, if we denote by $\alpha_n$ the image of
$\alpha$ in $U_1/U_n$, then $(\alpha_n)^{p^{n-2}}\neq1$ whereas
$(\alpha_n)^{p^{n-1}}=1$, or equivalently $\alpha_n$ has order
$p^{n-1}$.  But the group $U_1/U_n$ has order $p^{n-1}$ because it has
a filtration
$$
U_n/U_n\subset U_{n-1}/U_n\subset\cdots\subset U_2/U_n\subset U_1/U_n
$$
whose successive quotients $(U_i/U_n)/(U_{i+1}/U_n)=U_i/U_{i+1}$ have
order~$p$ (??), so $U_1/U_n$ is cyclic (and $\alpha_n$ is a
generator).

The case $p=2$ is similar.  Choose an $\alpha$ (for example
$\alpha=1+2^2$) such that $\alpha\in U_2$ but $\alpha\notin U_3$.  We
observe (??) that $\alpha^{2^i}\in U_{i+2}$ but $\alpha^{2^i}\notin
U_{i+3}$ for every $i>0$, so if $\alpha_n$ is the image of $\alpha$ in
$U_2/U_n$, then $(\alpha_n)^{2^{n-3}}\neq1$ whereas
$(\alpha_n)^{2^{n-2}}=1$.  But the group $U_2/U_n$ has order
$2^{n-2}$, so it is cyclic (and $\alpha_n$ is a generator). \cqfd

\numeroter {\it For\/ $p\neq2$, the group\/ $U_1$ is isomorphic to\/
  $\Z_p$.  For\/ $p=2$, the group\/ $U_2$ is isomorphic to\/ $\Z_2$,
  and\/ $U_1=\Z^\times\times U_2$.  More precisely, for every\/
  $\alpha\in U_1$ such that\/ $\alpha\notin U_2$ (resp.~$\alpha\in
  U_2$ such that\/ $\alpha\notin U_3$), there is a unique
  isomorphism\/ $f:\Z_p\to U_1$ (resp.~$\Z_2\to U_2$) such that\/
  $f(1)=\alpha$. }

{\it Proof}.  We have seen that the group $U_1/U_n$ (resp.~$U_2/U_n$)
is cyclic of order $p^{n-1}$ (resp.~$2^{n-2}$), and the image
$\alpha_n$, for any $\alpha$ such that $\alpha\in U_1$
(resp.~$\alpha\in U_2$) but $\alpha\notin U_2$ (resp.~$\alpha\notin
U_3$), is a generator.  Moreover, the image of $\alpha_{n+1}\in
U_1/U_{n+1}$ in $U_1/U_n$ is $\alpha_n$, so we get commutative
diagrams 
$$
\diagram{
\Z/p^n\Z&\droite{\sim}&U_1/U_{n+1}&\quad
\Z/2^{n-1}\Z&\droite{\sim}&U_2/U_{n+1}\cr
\vfl{\varphi_{n-1}}{}{5mm}&&\vfl{}{}{5mm}&\quad
\vfl{\varphi_{n-2}}{}{5mm}&&\vfl{}{}{5mm}\cr
\Z/p^{n-1}\Z&\droite{\sim}&U_1/U_n&\quad
\Z/2^{n-2}\Z&\droite{\sim}&U_2/U_n\cr
} 
$$
in which the horizontal maps are isomorphisms (coming from our choice
of $\alpha$).  Since $\Z_p$ is the projective limit of the vertical
maps on the left, and $U_1$ (resp.~$U_2$) is the projective limit of
the vertical maps on the right, it follows that there is a unique
isomorphism $\Z_p\to U_1$ such that $1\mapsto\alpha$.

Finally, for $p=2$, we have $U_1/U_2=(\Z/4\Z)^\times$ and the
isomorphism $\Z^\times\to(\Z/4\Z)^\times$, so the multiplication map
$\Z^\times\times U_2\to U_1$ is also an isomorphism~: clearly every
$x\in U_1$ can be uniquely written as $x=su$ ($s\in\Z^\times$, $u\in
U_2$).  \cqfd

\numeroter For $p=2$, the restriction of
$\lambda_8:\Z_2^\times\to\Z^\times$ to $U_2$ induces the unique
isomorphism $U_2/U_3\to\Z^\times$.  Conversely, $\lambda_8$ can be
recovered from this isomorphism via the projection
$\Z_2^\times\!/\Z^\times\to U_2$ (which can be written
$a\mapsto\lambda_4(a)a$).

\numeroter Note that the procedure in (??) gives each $U_m$ ($m>0$)
the structure of a (multiplicatively written) $\Z_p$-module.  Let us
explain this for $m=1$ for simplicity.  Let $a\in\Z_p$ and $u\in
U_1$~; they give rise to coherent sequences of elements
$a_n\in\Z/p^n\Z$ and $u_{n+1}\in U_1/U_{n+1}$, and hence a coherent
sequence of elements $u_{n+1}^{a_n}\in U_1/U_{n+1}$ which defines
$u^a\in U_1$ in the limit.  Similarly, for $m>1$, the system
$(U_m/U_{m+r+1}\to U_m/U_{m+r})_{r>0}$ consists of $\Z_p$-modules,
hence its projective limit $U_m$ is a $\Z_p$-module.  The $U_m$ are
all free of rank $1$ except when $p=2$ and $m=1$, in which case the
torsion subgroup is $\Z^\times\subset U_1$ (??).


\numeroter {\it Let\/ $p\neq2$ be a prime and let\/ $n>0$.  The
projection $G_{p^n}\to\F_p^\times$ has a canonical section, and\/
$G_{p^n}$ is canonically isomorphic to\/ $\F_p^\times\times(U_1/U_n)$.
The group\/ $U_1/U_n$ is cyclic of order~$p^{n-1}$ and it is generated
by\/ $1+p$.}

{\it Proof}.  Indeed, $G_{p^n}=\Z_p^\times/U_n$ (??), we know the
structure of $\Z_p^\times$ from (??), and we've seen that $U_1/U_n$ is
cyclic of order $p^{n-1}$ and generated by $1+p$ (??).  \cqfd

\numeroter {\it Let\/  $n>1$.  The projection
$G_{2^n}\to G_{2^2}$ has a canonical section, and\/ $G_{2^n}$ is
canonically isomorphic to\/ $\Z^\times\times U_2/U_n$. The group\/
$U_2/U_n$ is cyclic of order~$2^{n-2}$ and it is generated by\/
$1+2^2$.}

{\it Proof}.  Indeed, $G_{2^n}=\Z_2^\times/U_n$ (??), we know the
structure of $\Z_2^\times$ from (??), and we've seen that $U_2/U_n$ is
cyclic of order $2^{n-2}$ and generated by $1+2^2$ (??).  \cqfd

\numeroter It follows from the foregoing that $G_{p^n}$ ($n>0$) is
cyclic for every prime $p\neq2$, and $G_{2^n}$ ($n>1$) is cyclic if
and only if\/ $n=2$.

\numeroter The existence of the canonical section $\F_p^\times\to
G_{p^n}$ ($p$ odd prime, $n>0$) could also have been deduced from the
following algebraic lemma.  {\it Let\/ $0\to A\to E\to B\to0$ be an
exact sequence of finite commutative groups such that the orders\/
$a$, $b$ of $A$, $B$ have no common prime factors.  Let\/ $B'$ be the
set of\/ $x\in E$ such that\/ $bx=0$.  Then\/ $E$ is the internal
direct sum of\/ $A$ and\/ $B'$, and\/ $B'$ is the only subgroup of\/ $E$
isomorphic to\/ $B$.}

{\it Proof}.  Let $r,s\in\Z$ be such that $ar+bs=1$.  If $x\in A\cap
B'$, then $ax=bx=0$, hence $(ar+bs)x=x=0$, and cosequently $A\cap
B'=0$.  Next, since $bB=0$, we have $bE\subset A$, and hence $bsx\in
A$ for every $x\in E$.  At the same time, from $abE=0$ it follows that
$arx\in B'$ for every $x\in E$.  But every $x\in E$ can be written as
$x=arx+bsx$, so $E=A+B'$, and, in view of $A\cap B'=0$, the sum is
direct and the projection $E\to B$ induces an isomorphism $B'\to B$.
Conversely, if $B''$ is a subgroup of $E$ isomorphic to $B$, then
$bB''=0$, hence $B''\subset B'$ and $B''=B'$ because they have the
same order.  \cqfd

\numeroter  {\it For every\/ $a\in\Z$ prime to\/ $p$ and for every\/
$n>0$, the map\/ $(\ )^a:U_n\to U_n$ is an isomorphism.}

{\it Proof}.  Indeed, each $U_n$ is a $\Z_p$-module (??) and $a$,
being prime to $p$, is invertible in $\Z_p$ (??). \cqfd

\medbreak
\leftline{\it  Exercises}  
\medskip

\numeroter For which $u\in\Z_3^\times$ does the polynomial $T^3-u$
have a root in $\Z_3$~?  (The only cubes in $(\Z/3^2\Z)^\times$ are
$\pm1$, and the only cubes in $(\Z/3^3\Z)^\times$ are
$\pm1, \pm8, \pm10$. If $u$ is congruent to one of these $\pmod{3^3}$,
then (??) can be applied with $m=3$, $\delta=1$.)

\numeroter Let $a\in\F_p$, and let $x\in\Z_p$ be any lift of $a$ in
the sense that $\pi_1(x)=a$.  Show that the sequence
$(x^{p^n})_{n\in\N}$ converges in $\Z_p$ to $\omega(a)$.  ({\it Hint
  :\/} Reduce to the case $a=1$ and use the fact that $(\ )^p$ maps
$U_m$ into $U_{m+1}$).

\numeroter Let $R\subset\Z_p$ be the set of multiplicative
representatives (??) of $\F_p$.  Show that for every $(a_n)_{n\in\N}$
in $R^{\N}$, the series $\sum_{n\in\N}a_np^n$ converges in $\Z_p$ and
that the map $(a_n)_{n\in\N}\mapsto\sum_{n\in\N}a_np^n$ is a bijection
$R^{\N}\to\Z_p$.  ({\it Hint~:\/} Construct a reciprocal map by
associating to $x\in\Z_p$ the sequence $(a_n)_{n\in\N}$ inductively
defined by $a_0\equiv x\pmod p$, $a_1\equiv{(x-a_0)\over p}\pmod p$,
and so on.)

\numeroter Let $x,y\in\Z_p$, and let $(a_n)_{n\in\N}$,
$(b_n)_{n\in\N}$ be the sequences in $R^{\N}$ such that
$x=\sum_{n\in\N}a_np^n$ and $y=\sum_{n\in\N}b_np^n$.  Find the
sequences corresponding to $x+y$ and $xy$.

\numeroter  Let $\bar\F_p$ be an algebraic closure of the field $\F_p$, and give
the group $\bar\F_p^\times$ the discrete topology.  Show that every
continuous morphism of groups $\chi:\Z_p^\times\to\bar\F_p^\times$
factors through the quotient $\Z_p^\times\to\F_p^\times$.  ({\it
Hint~:\/} As $\chi$ is continuous and $\bar\F_p^\times$ is discrete,
there is an $n>0$ such that $U_n\subset\Ker\chi$.  Note that every
element of $\bar\F_p^\times$ has order prime to~$p$.)

\numeroter Let $G$ be a finite commutative group (written
multiplicatively), and let $s$ be the product of all elements of $G$.
Show that if $G$ has a unique element $\tau$ of order~$2$, then
$s=\tau$, otherwise $s=1$.  (Consider the involution $x\mapsto x^{-1}$
of $G$.)  Take $G=G_m$ for some integer $m>0$.  Show that $s=-1$ if
$m=4$ or $m=l^a$ or $m=2l^a$ for some prime~$l\neq2$~; otherwise
$s=1$.  (Use the structure of $G_{p^n}$ (??) and the Chinese remainder
theorem.)

\numeroter Recall that for every integer $m>0$ we are using the
abbreviation $G_m=(\Z/m\Z)^\times$.  A {\it character\/} of $G_m$ is a
morphism $\chi:G_m\to\C^\times$ of groups~; quadratic characters (??)
are a special case.  If $m$ is a multiple of some $m'>0$, then every
character of $G_{m'}$ can be viewed as a character of $G_m$ via the
projection $\varphi_{m,m'}:G_m\to G_{m'}$.  If $\chi$ is {\it not\/}
of the form $\chi'\circ\varphi_{m,m'}$ for any character $\chi'$ of
$G_{m'}$ for any divisor $m'<m$, we say that $\chi$ is {\it
primitive}, and that the {\it conductor\/} of $\chi$ is $m$.  The
quadratic characters $\lambda_p$ ($p$ odd prime), $\lambda_4$,
$\lambda_8$, $\lambda_4\lambda_8$ are primitive of conductors $p$,
$4$, $8$, $8$ respectively~; show that these are the only primitive
quadratic characters of prime-power conductor (???).

\numeroter For every character $\chi$ of $G_m$, we have $\chi(-1)^2=1$
and hence $\chi(-1)=1$ or $\chi(-1)=-1$.  We say that $\chi$ is {\it
even\/} if $\chi(-1)=1$, {\it odd\/} if $\chi(-1)=-1$~; this extends
the terminology introduced earlier (???)  for quadratic characters.
Define $\varepsilon(\chi)\in\Z_2$ by
$\chi(-1)=(-1)^{\varepsilon(\chi)}$, so that $\varepsilon(\chi)=0$ if
$\chi$ is even and $\varepsilon(\chi)=1$ if $\chi$ is odd.  (For the
quadratic characters $\lambda_*$, we had earlier defined
$\varepsilon_*(-1)\in\F_2$ (??)~; the newly defined
$\varepsilon(\lambda_*)$ is the multiplicative representative (??) of
$\varepsilon_*(-1)$.  The advantage of the new definition is that
whereas the congruence
$\varepsilon(\chi_1\chi_2)-\varepsilon(\chi_1)-\varepsilon(\chi_2)\equiv0\pmod2$
holds for any two characters $\chi_1,\chi_2$, we have the inequality
$\varepsilon(\chi_1\chi_2)-\varepsilon(\chi_1)-\varepsilon(\chi_2)\neq0$
when both $\chi_1$ and $\chi_2$ are odd.) Let $\chi_1,\chi_2$ be
primitive quadratic characters of $G_{m_1}$, $G_{m_2}$ respectively
such that $\gcd(m_1,m_2)=1$, and let $\chi=\chi_1\chi_2$ be the
quadratic character of $G_{m_1m_2}$ coming from the isomorphism
$G_{m_1m_2}\to G_{m_1}\times G_{m_2}$~; show that $\chi$ is primitive.
Show that for every primitive $4$-th root $i$ of $1$ in $\C$, we have
$$
\chi_1(m_2)\chi_2(m_1)=i^{\varepsilon(\chi_1\chi_2)-\varepsilon(\chi_1)-\varepsilon(\chi_2)}. 
$$

\numeroter  We have seen (???) that the identity $(T+1)^p=T^p+1$ holds in
the polynomial ring $\F_p[T]$.  Show that, conversely, if
$(T+1)^n=T^n+1$ in $\Z/n\Z[T]$ for some $n>1$, then $n$ is prime.  (If
$n$ is not prime, let $p<n$ be a prime divisor of $n$ and put
$\delta=v_p(n)$, so that $\delta>0$.  The coefficient of $T^p$ in the
binomial expansion of $(T+1)^n$ is ${n\choose
p}=n(n-1)\ldots(n-(p-1))/p!$, so $v_p({n\choose p})=\delta-1$.  This
implies that ${n\choose p}\not\equiv0\pmod n$.)

\vfill\eject

\centerline{\bf Lecture 5}
\medskip
\centerline{$\Q_p$}
\bigskip

\numeroter Let $p$ be a prime number.  We define $\Q_p$ to be the
field of fractions of $\Z_p$ (??).  Since every $x\neq0$ in $\Z_p$ can
be uniquely written as $x=p^mu$ ($m\in\N$, $u\in\Z_p^\times$), we have
$\Q_p=\Z_p[{1\over p}]$, and every $x\neq0$ in $\Q_p$ can be uniquely
written as $x=p^mu$ ($m\in\Z$, $u\in\Z_p^\times$)~; we pose
$v_p(x)=m$.  This new definition of $v_p$ extends our earlier
definition (??).

\numeroter The resulting homomorphism $v_p:\Q_p^\times\to\Z$ is a {\it
  valuation\/} in the sense that $v_p(x+y)\ge\Inf(v_p(x), v_p(y))$ for
all $x,y\in\Q_p$ (with the convention $v_p(0)=+\infty$).  For
$x\in\Q_p^\times$, put $|x|_p=p^{-v_p(x)}$.  Then $|\ |_p$ is a
homomorphism $\Q_p^\times\to\R^\times_+$ satisfying
$|a+b|_p\le\Sup(|a|_p+|b|_p)$ (with the convention $|0|_p=0$), and
$d_p(x,y)=|x-y|_p$ is a distance on $\Q_p$ satisfying
$d_p(x,z)\le\Sup(d_p(x,y),d_p(y,z))$, for all $x,y,z\in\Q_p$.

\numeroter For every $x\in\Q^\times$, we
have $|x|_v=1$ for almost all places $v$ of $\Q$ and the {\it product
  formula\/} $\prod_v|x|_v=1$ holds.  Indeed, in view of the
  multiplicativity of the $|\ |_v$, it is enough to verify it for
  $x=-1$ and $x=p$ for every prime $p$.

\numeroter For every $m\in\Z$, we have the sub-$\Z_p$-module of $\Q_p$
generated by $p^m$, and the inclusion $p^{m+1}\Z_p\subset p^m\Z_p$~;
the union of this increasing sequence (when $m\to-\infty$) is $\Q_p$.

\numeroter {\it The field\/ $\Q_p$ is locally compact, complete for
  $d_p$, and the subring\/ $\Z[{1\over p}]$ is dense.}

{\it Proof}.  Recall (??) that $\Z_p$ is compact~; as it is defined as
a subspace of $\Q_p$ by $d_p(0,x)<p$, it is an open neighbourhood of
$0$.  Therefore $\Q_p$ is locally compact and therefore complete, like
any locally compact commutative group.  Another way to prove
completeness is to remark that if $(x_n)$ is a fundamental sequence in
$\Q_p$, then there is an $M>0$ such that the sequence $(p^Mx_n)$ is
(fundamental and) in $\Z_p$~; if $y\in\Z_p$ is the limit of the latter
sequence, then the former sequence has the limit $x=p^{-M}y$ in
$\Q_p$.

\numeroter We have seen (??) that every $x\in\Z_p$ can be uniquely
written as $\sum_{n\in\N}b_np^n$, with $b_i\in[0,p[$.  It follows that
    every $x\in\Q_p$ can be uniquely written as $\sum_{n\ge
      v_p(x)}b_np^n$, making $x$ the limit of the sequence $(s_m)_m$
    of partial sums of the series representing $x$.  But each $s_m$ is
    in $\Z[{1\over p}]$, so $\Z[{1\over p}]$ is dense in $\Q_p$.
    \cqfd

\numeroter An element $\pi\in\Q_p$ is called a {\it uniformiser\/} if
$v_p(\pi)=1$~; the simplest example is $\pi=p$.  The choice of a
uniformiser $\pi$ leads to the splitting $1\mapsto\pi$ of the short
exact sequence
$$
1\to\Z_p^\times\to\Q_p^\times\to\Z\to0,
$$
and thus to an isomorphism $(m,u)\mapsto\pi^mu$ of groups
$\Z\times\Z_p^\times\to\Q_p^\times$.  If we also choose a generator
$\alpha$ of the (free rank-$1$) $\Z_p$-module $U_1$ when $p\neq2$
(resp.~$U_2$ when $p=2$), then
$\Q_p^\times=\Z\times\F_p^\times\times\Z_p$ ($p\neq2$) (???), and
$\Q_2^\times=\Z\times\Z^\times\times\Z_2$ (???).

\numeroter {\it If \/ $p\neq2$, the group\/ $\Q_p^\times/\Q_p^{\times
    2}$ consists of $\{\bar 1;\bar u; \bar p,\bar u\bar p\}$, where
  $u\in\Z_p^\times$ is any unit such that\/ $\lambda_p(u)=-1$.}

{\it Proof}.  This follows immediately from the 
  isomorphism (??) and the fact that $\Z_p^\times/\Z_p^{\times 2}=\{1,\bar
  u\}$ (??).

\numeroter {\it We have\/  $\Q_2^\times/\Q_2^{\times 2}=
\{\bar 1;\bar 5; -\bar 1,-\bar 5; \bar 2, \overline{10},
-\bar2, -\overline{10}\}$.}

{\it Proof}.  This follows similarly from the isomorphism (??) and the
fact that $\Z_2^\times/\Z_2^{\times 2}=\{\bar 1;\bar 5; -\bar 1,-\bar
5\}$ (??).  \cqfd

\numeroter For every $p$, the morphism $\nu_p(x)=(-1)^{v_p(x)}$ is a
quadratic character of $\Q_p^\times$~; we call it the {\it
  unramified\/} quadratic character.  Choosing a uniformiser $\pi$, we
get a retraction $x\mapsto x\pi^{-v_p(x)}$ of the inclusion
$\Z_p^\times\to\Q_p^\times$, allowing us to view quadratic characters
of $\Z_p^\times$ as {\it ramified\/} quadratic characters of
$\Q_p^\times$~; for the moment, don't worry about the meaning of these
words.  We choose the uniformiser $p$ to fix ideas.

\numeroter {\it For $p\neq2$, the quadratic characters\/ $\nu_p$,
  $\lambda_p$ constitute a basis of the\/ $\F_2$-space
  $\Hom(\Q_p^\times,\Z^\times)$~; their values on the basis\/ $\bar
  u,\bar p$ of\/ $\Q_p^\times/\Q_p^{\times2}$ are given by the matrix}
$$\bordermatrix{&\nu_p&\lambda_p\cr
u&\phantom{+}1&-1\cr
p&-1&\phantom{+}1\cr
}.$$  \cqfd

\numeroter {\it
The quadratic characters\/ $\nu_2, \lambda_4,\lambda_8$ constitute a
basis of the\/ $\F_2$-space \/ $\Hom(\Q_2^\times,\Z^\times)$.
Their values on the basis\/ $\bar5,-\bar1,\bar2$ of\/
$\Q_2^\times/\Q_2^{\times2}$ are given by }
$$\bordermatrix{&\nu_2&\lambda_4&\lambda_8\cr
\phantom{-}5&\phantom{+}1&\phantom{+}1&-1\cr
-1&\phantom{+}1&-1&\phantom{+}1\cr
\phantom{-}2&-1&\phantom{+}1&\phantom{+}1\cr
}.
$$ \cqfd

\numeroter Let us extend the hilbertian symbol $(\ ,\ )_p$ to $\Q_p$.
Let $a,b\in\Q_p^\times$.  Write
$a=p^{v_p(a)}u_a$ and $a=p^{v_p(b)}u_b$, and as in (??) put
$$
t_{a,b}
=(-1)^{v_p(a)v_p(b)}u_a^{v_p(b)}u_b^{-v_p(a)}
$$
If $p\neq2$, define
$
(a,b)_p=\lambda_p(t_{a,b})
$, 
as in (??).  Similarly, as in (??),  define 
$$
(a,b)_2
=(-1)^{\varepsilon_4(u_a)\varepsilon_4(u_b)+\varepsilon_8(t_{a,b})}
$$

\numeroter It is readily verified this new definition is compatible
with the old (??), and continues to enjoy all the properties listed in
(??) to (??).

\numeroter {\it Let\/ $a,b\in\Q_p^\times$.  If\/ $a+b=1$, then
  $(a,b)_p=1$.}

{\it Proof}.  ??? \cqfd

\numeroter For every $p$, the pairing $(\ ,\ )_p$ is {\it
  invertible\/} in the sense that its matrix (for $p\neq2$ and $p=2$
  respectively)
$$\def\\{\phantom{+}}
\bordermatrix{&u&p\cr u&\\1&-1\cr p&-1&\lambda_4(p)\cr},\qquad
\bordermatrix{&\\5&-1&\\2\cr
\\5&\\1&\\1&-1\cr
-1&\\1&-1&\\1\cr
\\2&-1&\\1&\\1\cr},
$$
(with entries in $\Z^\times$) with respect to the given basis of
$\Q_p^\times/\Q_p^{\times 2}$ is invertible (when viewed with entries
in the field $\F_2$).  Indeed, in $\F_2$,
$$
\qquad
\left|\matrix{0&1\cr 1&\varepsilon_4(p)\cr}\right|=1,
\qquad\qquad\quad
\left|\matrix{0&0&1\cr 0&1&0\cr 1&0&0\cr}\right|=1.
$$

\numeroter The reader must have noticed that these matrices are the
same as the ones in (??)  and~(??) giving the values of the basic
quadratic characters on the chosen basis of
$\Q_p^\times/\Q_p^{\times2}$, with the important exception of the
entry $\lambda_4(p)=\lambda_p(-1)$.  This phenomenon will get
explained later when we will have interpreted the hilbertian symbol in
terms of the reciprocity isomorphism for the maximal abelian extension
of $\Q_p$ of exponent~$2$.

\numeroter Before proceeding further, we need a small general lemma.
Let $k$ any field of characteristic $\neq2$, and let $a,b\in
k^\times$.  If there is a pair $(x,y)\in k^2$ such that $ax^2+by^2=1$,
then certainly there is a triple $(x,y,s)\neq(0,0,0)$ such that
$ax^2+by^2=s^2$.  {\it Conversely, if there is such a triple, then
  there is a desired pair.}  This is clear if $s\neq0$.  If $s=0$ (in
which case $x\neq0$ and $y\neq0$), we have $-a=b/t^2$ with $t=y/x$ and
$a\left({a+1\over2a}\right)^2 + b\left({a-1\over2at}\right)^2=1$, so a
suitable pair exists.

\numeroter {\it Let\/ $a,b\in\Q_p^\times$.  There exist\/ $x,y\in\Q_p$
  such that\/ $ax^2+by^2=1$ if and only if\/ $(a,b)_p=1$.}

{\it Proof}.  Suppose that there do exist $x,y\in\Q_p$ such that
$ax^2+by^2=1$.  If $x=0$ (resp.~$y=0$), then $b$ (resp.~$a$) is in
$\Q_p^{\times2}$, and hence $(a,b)_p=1$.  If $xy\neq0$, then
$$
(a,b)_p=(ax^2,by^2)_p=1,
$$
since $(c,d)_p=1$ whenever $c+d=1$ (???).  It remains to prove the
converse. 

\numeroter Suppose that $(a,b)_p=1$~; we have to show that there do
exist $x,y\in\Q_p$ such that $ax^2+by^2=1$.  Since the value $(a,b)_p$
as well as the existence of $x,y$ depend only on the classes of $a$
and $b$ modulo $\Q_p^{\times2}$, we need only consider the following
cases
$$
v_p(a)=v_p(b)=0\;;\quad
v_p(a)=0,\ v_p(b)=1\;;\quad
v_p(a)=v_p(b)=1.
$$
\numeroter In fact, the last case can be reduced to the middle case upon
replacing $a$ by $-ab^{-1}$.  First,
$$
(-ab^{-1},b)_p=(a,b)_p(-b^{-1},b)_p=(a,b)_p(-b,b)_p=(a,b)_p.
$$
Secondly, the existence of $x,y\in\Q_p$ such that
$-ab^{-1}x^2+by^2=1$ is equivalent to the existence of
$(x,y,z)\neq(0,0,0)$ in $\Q_p^3$ such that $-ab^{-1}x^2+by^2=z^2$ (??).
Multiplying throughout by $b$ and rearranging, the latter becomes
$ax^2+bz^2=(by)^2$, which can be seen as before to be equivalent to
the existence of $x,y\in\Q_p$ such that $ax^2+by^2=1$.

\numeroter {\it The case\/} $v_p(a)=v_p(b)=0$~; {\it subcase\/}
$p\neq2$.  Consider the subset $S\subset\F_p$ of all elements of the
form $\bar a\xi^2$ ($\xi\in\F_p$) and the subset $T\subset\F_p$ of all
elements of the form $1-\bar b\eta^2$ ($\eta\in\F_p$)~; each of these has
$(p+1)/2$ elements, so there exist $x,y\in\Z_p$ such that
$ax^2+by^2\equiv1\pmod p$. If $x\not\equiv0\pmod p$, then the unit
$(1-by^2)a^{-1}$ is a square $\pmod p$, and hence the square of some
$t\in\Z_p^\times$.  We then have $at^2+by^2=1$, and we are done.  If
$x\equiv0\pmod p$, then for the same reason $b=t^2$ for some
$t\in\Z_p^\times$, and we have $a.0^2+b.(t^{-1})^2=1$ in $\Q_p$.

\numeroter {\it The subcase\/} $p=2$.  As
$(a,b)_2=(-1)^{\varepsilon_4(a)\varepsilon_4(b)}=1$, we may suppose
(up to interchanging $a$ and $b$) that $a\equiv1\pmod4$, or
equivalently $a\equiv_8 1$ or $a\equiv_8 5$.  If $a\equiv1\pmod8$, then
there is a $t\in\Z_2^\times$ such that $a=t^2$, and hence
$a(t^{-1})^2+b.0^2=1$ in $\Q_2$.  Suppose finally that
$a\equiv5\pmod8$.  As $4b\equiv4\pmod8$, we have $a\equiv 1-4b\pmod8$,
and there is a $t\in\Z_2^\times$ such that $t^2=(1-4b)a^{-1}$.  We
then have $a.t^2+b.2^2=1$, and we are done.

\numeroter {\it The case\/} $v_p(a)=0$, $v_p(b)=1$~; {\it subcase\/}
$p\neq2$.  The hypothesis $(a,b)_p=\lambda_p(a)=1$ implies that $\bar
a\in\F_p^{\times 2}$, so (??) there is a $t\in\Z_p^\times$ such that
$a=t^2$, and then $a(t^{-1})^2+b.0^2=1$.

\numeroter {\it The subcase\/} $p=2$.  The reader should check that
the hypothesis $(a,b)_2=1$ is equivalent in this subcase to
``\thinspace either $a\equiv1\pmod8$ or
$a\equiv1-b\pmod8$\thinspace'', so there exists a $t\in\Z_2^\times$
such that $t^2=a$ or $t^2=(1-b)a^{-1}$.  In the former case we have
$a(t^{-1})^2+b.0^2=1$, and in the latter $at^2+b.1^2=1$, so we are
done. \cqfd

\numeroter Let $b\in\Q_p^\times$, and put $K_b=\Q_p(\sqrt b)$, so that
the degree $[K_b:\Q_p]$ equals $1$ or $2$ according as
$b\in\Q_p^{\times 2}$ or $b\notin\Q_p^{\times 2}$.  We have the norm
homomorphism $N_b:K_b^\times\to\Q_p^\times$ which is the identity in
case $[K_b:\Q_p]=1$ and sends $x+y\sqrt b$ ($x,y\in\Q_p$) to
$x^2-by^2$ in case $[K_b:\Q_p]=2$.

\numeroter {\it Let\/ $a,b\in\Q_p^\times$.  We have\/ $(a,b)_p=1$ if
  and only if\/ $a\in N_b(K_b^\times)$.}

{\it Proof}.  There is nothing to prove if $b\in\Q_p^{\times2}$.  If
not, the proposition follows from the equivalence of the following
four conditions~: \par\noindent 
{\it i\/}) $(a,b)_p=1$,\hfill\break 
{\it ii\/}) there exists a pair $(x,y)$ in $\Q_p^2$ such that
  $ax^2+by^2=1$ (??),\hfill\break 
{\it iii\/}) there exists a triple $(x,y,z)\neq(0,0,0)$ in $\Q_p^3$ such that
  $ax^2+by^2=z^2$,\hfill\break 
{\it iv\/}) there exists a pair $(y,z)$ in $\Q_p^2$ such that
$a=z^2-by^2$ (??). \cqfd

\numeroter This proposition can be interpreted as saying that
$N_b(K_b^\times)=b^\perp$, where the orthogonal is taken with respect
to the hilbertian pairing.

\numeroter Note finally that for $a,b\in\R^\times$ and $K_b=\R(\sqrt
b)$, we have $a\in N_b(K_b^\times)$ if and only if $(a,b)_\infty=1$.

\numeroter {\it For every quadratic extension\/ $E$ of\/ $\Q_p$, the
  subgroup\/ $N_{E|\Q_p}(E^\times)$ of \/$\Q_p^\times$ is an open
  subgroup of index~$2$.} \cqfd

\bigbreak 
\leftline{\it Exercises}
\medskip

\numeroter Show that if $v_p(x)<v_p(y)$, then
$v_p(x+y)=v_p(x)$.  Give an example where $v_p(x)=v_p(y)$ but
$v_p(x+y)>v_p(x)$.

\numeroter Show that for every $n\in\Z$, the subset
$v_p^{-1}([n,+\infty])$ of $\Q_p$ is $p^n\Z_p$.

\numeroter  The product formula (???) can be rewritten as
$|x|_\infty=\prod_p|x|_p^{-1}$ ($x\in\Q^\times$).  Show the unicity of
this formula in the following sense~: if the reals $s_p>0$ are such
that $|x|_\infty=\prod_p|x|_p^{-s_p}$ for every $x\in\Q^\times$, then
$s_p=1$ for every $p$.  (Take $x=l$, where $l$ runs through the
primes).

\numeroter Show that the polynomial $(T^2-2)(T^2-17)(T^2-34)$ has 
roots in $\R$ and in $\Q_p$ for every prime $p$ but doesn't have a
root in $\Q$.

\numeroter For every integer $m>0$, find an $\F_2$-basis of
$\Hom(G_m,\Z^\times)$, where $G_m=(\Z/m\Z)^\times$. ({\it Hint :\/}
The structure of $G_{p^n}$ was determined in (??)~; note that
$\lambda_p$ is a quadratic character of $G_m$ for every odd prime
divisor $p$ of $m$, so is $\lambda_4$ if $v_2(m)>1$, and so is
$\lambda_8$ if $v_2(m)>2$.)


\numeroter Let $a,b,c\in\Z_2^\times$.  Show that a necessary and
sufficient condition for $ax^2+by^2+cz^2=0$ to have only the trivial
solution $(0,0,0)$ in $\Q_2$ is that $a\equiv b\equiv c\pmod4$.  ({\it
  Hint~:} The equation $dx^2+ey^2=1$ ($d,e\in\Q_2^\times$) has a
solution in $\Q_2$ if and only if $(d,e)_2=1$.  If
$d,e\in\Z_2^\times$, then $(d,e)_2=(-1)^{{d-1\over2}{e-1\over2}}$.
Finally take $d=-ca$, $e=-cb$ and show that the condition $(d,e)_2=-1$
is equivalent to the given condition $a\equiv b\equiv c\pmod4$.)

\numeroter Show that for $x\in\Q_p^\times$, the relation
``$x\in\Z_p^\times$'' (or ``$v_p(x)=0$'') is equivalent to
``$x\in\Q_p^{\times m}$ for every $m$ prime to $(p-1)p$''.

\numeroter Show that the only morphism $\sigma:\Q_p\to\Q_p$ of 
fields is the identity. (Conclude from (???) that
$\sigma(\Z_p^\times)\subset\Z_p^\times$, so $\sigma$ preserves the
valuation and $\sigma(p^r\Z_p)\subset p^r\Z_p$ for every $r\in\Z$, and
hence $\sigma$ is continuous.)

\numeroter Show that the only morphism $\sigma:\R\to\R$ of 
fields is the identity. (As $\R^{\times 2}=\R^\times_+$, the order is
preserved by $\sigma$ and hence it is continuous.)

\numeroter  For every prime $p$ and every $m>0$, compute the index of
the subgroup $\Q_p^{\times m}\subset\Q_p^\times$.  (Use the
decomposition (???).)

\numeroter Let $x\in\Q_p$, let $x=\sum_{n\in\Z}a_np^n$ be its
$p$-adic exapnsion (so that $a_n\in[1,p[$ and $a_n=0$ for almost all
$n<0$), put $\langle x\rangle_p=\sum_{n<0}a_np^n$, and view $\langle
x\rangle_p$ as a real number in $[0,1[\;\cap\;\Z[1/p]$.  Show that
$x\mapsto e^{2i\pi\langle x\rangle_p}$ is a continuous morphism of
groups $\Q_p\to\C^\times$ whose kernel is $\Z_p$.

\numeroter Show that the image of every {\it continuous\/}
morphism $\psi:\Q_p\to\C^\times$ (of locally compact groups) is
contained in the subgroup of $p$-power roots of~$1$.  (For every
$m\in\Z$, the subgroup $p^m\Z_p$ is compact, so its image
$\psi(p^m\Z_p)$ is contained in the unit circle $\U$.  As $\Q_p$ is
the union of these compact subgroups, $\psi(\Q_p)\subset\U$.  Next,
let $V$ be an open neighbourhood of $1\in\C^\times$ which contains
only the trivial subgroup of $\C^\times$.  Then $\psi^{-1}(V)$ is an
open neighbourhood of $0\in\Q_p$ and hence contains $p^N\Z_p$ for some
(and in fact for every) sufficiently large $N\in\Z$, so that
$\psi(p^N\Z_p)=1$~; fix such an $N$.  Finally observe that
$\Q_p/p^N\Z_p$ is the union of $p^m\Z_p/p^N\Z_p$ for $m<N$, and each
of these quotients is (cyclic) of order $p^{N-m}$.  Conclude.)

\numeroter Show that for every continous morphism of groups
$\chi:\Q_p^\times\to\C^\times$, there exists an $n>0$ such that
$\chi(U_n)=1$. We say that $\chi$ is {\it unramified\/} if
$\chi(\Z_p^\times)=1$, (at worst) {\it tamely ramified\/} if
$\chi(U_1)=1$, {\it totally ramified\/} if $\chi(\pi)=1$ for some
uniformiser $\pi$ of $\Q_p$, and {\it wildly ramified\/} if
$\chi(U_1)\neq1$.  Classify all unramified $\chi$, all tamely ramified
$\chi$, and all totally ramified $\chi$.  (Use the decomposition (???)
of $\Q_p^\times$.)

\numeroter (Bourgain-Larsen, 2014)  For every subgroup
$G\subset\Q^\times$ and every finte set $S$ of places of $\Q$, denote
by $G_S$ the closure of $G$ in the product $\prod_{v\in
S}\Q_v^\times$~; in particular, for every place $v$ of $\Q$, the
closure of $G$ in $\Q_v^\times$ is denoted $G_v$.  We want to give an
example of a subgroup $G$ (of finite index in $\Q^\times$) and a
triple $a,b,c\in\Q^\times$ such that the equation $ax+by+cy=0$ has a
solution $x_v,y_v,z_v$ in $G_v$ at every place $v$ of $\Q$ but no
solution $x,y,z$ in $G$.

\numeroter  Take $G\subset\Q^\times$ to be the subgroup generated by
the numbers $3^m5^n$ such that $m\equiv n\pmod4$ and all numbers
$t\in\Q^\times$ which, when viewed in $\Q_3^\times$ and $\Q_5^\times$,
are such that $t\in1+3\Z_3$ and $t\in1+5\Z_5$ (so that in particular
$v_3(t)=0$ and $v_5(t)=0$).  Show that the index of $G$ in $\Q^\times$
is $4.\varphi(3).\varphi(5)=32$.  Show that $G_v=\Q_v^\times$ for
every place $v\neq3,5$, and that $G_3=3^\Z(1+3\Z_3)$,
$G_5=5^\Z(1+5\Z_5)$.

\numeroter  However, $G_{\{3,5\}}$ is not equal to the product
$G_3\times G_5$~; rather, it is equal to the subgroup consisting of
$(t_3,t_5)\in G_3\times G_5$ such that $v_3(t_3)\equiv
v_5(t_5)\pmod4$.

\numeroter Show that the equation $63x+30y+25z=0$ has a solution 
in $G_3$ (for example $-5,3,9$), and that every solution $x_3,y_3,z_3$
satisfies
$$
v_3(x_3)=v_3(y_3)-1=v_3(z_3)-2.
$$
Similarly, show that the given equation has a solution in $G_5$ (for
example $25, -45, -9$), and that every solution $x_5,y_5,z_5$
satisfies
$$
v_5(x_5)=v_5(y_3)+1=v_5(z_5)+2.
$$

\numeroter  Conclude that the equation  $63x+30y+25z=0$ has no
solutions $x,y,z\in G$ (even though it has a solution in $G_v$ for
every place $v$ of $\Q$) because it doesn't have any solution in
$G_{\{3,5\}}$.

\vfill\eject

\centerline{\bf Lecture 6}
\medskip
\centerline{$\Q_p\big(\sqrt{\Q_p^\times}\big)$}
\bigskip

\numeroter Let $p$ be an {\it odd\/} prime number (resp.~$p=2$).  We
have seen that the field $\Q_p$ has three (resp.~seven) quadratic
extensions, namely those obtained by adjoining $\def\\#1{\sqrt{#1}}
\\u,\;\\{-p},\\{u.-p}$ (the reason for choosing $-p$ instead of $p$ is
that $(-p,p)_p=1$ whereas $(p,p)_p=1$ only when $\lambda_4(p)=1$)
$$\def\\#1{\sqrt{#1}}
 (\hbox{resp. }\\5,\;\\{-1},\\{5.-1},\;\\2,\;\\{5.2},\;\\{-1.2},\\{5.-1.2}).
$$
where, for $p\neq2$, any unit $u\in\Z_p^\times$ such that
 $\lambda_p(u)=-1$ can be chosen.  We have also seen that the group
 $\Q_p^\times$ has three (resp.~seven) quadratic characters, namely
$\nu_p,\;\lambda_p,\;\nu_p\lambda_p$
$$
(\hbox{resp. }
\nu_2,\;\;\lambda_4,\;\;\nu_2\lambda_4,\;\;\lambda_8,\;\;
\nu_2\lambda_8,\;\;\lambda_4\lambda_8,\;\;\nu_2\lambda_4\lambda_8)
$$
(with the choice of $p$ as a uniformiser of $\Q_p$).  One gets the
feeling that there is a canonical bijection between these two lists,
but what characterises this bijection~?

\numeroter {\it There is a unique bijection\/ $E\mapsto\chi$ between
  the set of quadratic extensions of\/ $\Q_p$ and the set of quadratic
  characters of\/ $\Q_p^\times$ such that\/
  $\Ker(\chi)=N_{E|\Q_p}(E^\times)$.  The induced map\/
  $\Q_p^\times\!/\Q_p^{\times2}\to\Hom(\Q_p^\times,\Z^\times)$ is an
  isomorphism of groups.}

{\it Proof}.  We claim that the bijection given by the order in which
the two lists have been written in (??) is the required one.  This is
simple to verify~; let us check for example that $E=\Q_p(\sqrt u)$
(resp.~$E=\Q_2(\sqrt5)$) corresponds to $\chi=\nu_p$.  As
$\Ker(\nu_p)=p^{2\Z}.\Z_p^\times$, this amounts to checking that an
$a\in\Q_p^\times$ is in $N_{E|\Q_p}(E^\times)$ if and only if
$v_p(a)\in2\Z$.  Now, $a\in N_{E|\Q_p}(E^\times)$ if and only if there
exist $x,y\in\Q_p$ such that $a=x^2-uy^2$ (resp.~$a=x^2-5y^2$), which
can be seen to be equivalent to $(a,u)_p=1$ (resp.~$(a,5)_2=1$).  The
explicit formul\ae\ for the pairing $(\ ,\ )_p$ tell us that this last
condition holds if and only if $v_p(a)\in2\Z$.

The proof in the other cases is similar.  For example, to check that
$\chi=\lambda_p$ corresponds to $E=\Q_p(\sqrt{-p})$ (for $p\neq2$), we
have to show that for a given $a=p^mu$ (with $m\in\Z$ and
$u\in\Z_p^\times$), we have $a\in N_{E|\Q_p}(E^\times)$ if and only if
$\lambda_p(u)=1$. The first condition translates into the existence of
$x,y\in\Q_p$ such that $a=x^2+py^2$, or equivalently $(a,-p)_p=1$.
But $(a,-p)_p=(p^mu,-p)_p=(p^m,-p)_p(u,-p)_p=\lambda_p(u)$, so we are
done.  

As a random example, take $\chi=\lambda_8$ and $E=\Q_2(\sqrt2)$. An
$a\in\Q_2^\times$ is in $N_{E|\Q_2}(E^\times)$ if and only if there
exist $x,y\in\Q_2$ such that $a=x^2-2y^2$, which is equivalent to
$(a,2)_2=1$, which is equivalent to $\lambda_8(a)=1$.  The
conscientious reader should check the remaining cases. 
 
We thus get a bijection
$\Q_p^\times\!/\Q_p^{\times2}\to\Hom(\Q_p^\times,\Z^\times)$, and it
remains to show that it is an isomorphism.  A cursory look at the two
lists is enough to conclude that it is indeed so.  \cqfd

\numeroter Actually, (??) is a corollary of what we have proved
earlier.  We have canonical bijections between the following four
sets~:

$\cal E$, the set of quadratic extensions of $\Q_p$,

$\cal H$, the set of index-$2$ subgroups of 
  $\Q_p^\times\!/\Q_p^{\times 2}$,

$\cal L$, the set of order-$2$ subgroups of 
  $\Q_p^\times\!/\Q_p^{\times 2}$,

$\cal Q$, the set of quadratic characters of $\Q_p^\times$.

The bijection ${\cal E}\to{\cal H}$ sends $E$ to
$N_{E|Q_p}(E^\times)$, the bijection ${\cal H}\to{\cal L}$ sends $H$
to $H^\perp$, where the perpendicular is taken with respect to the
hilbertian pairing $(\ ,\ )_p$, the bijection ${\cal L}\to{\cal E}$
sends $L$ to $\Q_p(\sqrt L)$, and the bijection ${\cal Q}\to{\cal H}$
sends $\chi$ to $\Ker(\chi)/\Q_p^{\times2}$.

\numeroter The unramfied quadratic character $\nu_p$ corresponds to
$\Q_p(\sqrt u)$ (for $p\neq2$) and to $\Q_2(\sqrt5)$ (for $p=2$).
These quadratic extensions will therefore be called {\it unramified\/}
(over $\Q_p$).  Note that $\Q_2(\sqrt 5)$ contains $\sqrt{-3}$ and
hence a primitive $3$-rd root of~$1$.  Similarly, we shall see later
for $p\neq2$ that $\Q_p(\sqrt u)$ contains a primitive $(p^2-1)$-th
root of~$1$.

\numeroter Let $M$ be the maximal abelian extension of exponent~$2$
(or equivalently the compositum of all quadratic extensions) of
$\Q_p$.  Concretely, $M=\Q_p(\sqrt u,\sqrt{-p})$ if $p\neq2$ and
$M=\Q_2(\sqrt5,\sqrt{-1},\sqrt{2})$ for $p=2$.  Let $G=\Gal(M|\Q_p)$.
As an $\F_2$-space, the dimension of $G$ is $2$ (resp.~$3$), so there
are many isomorphisms of $G$ with the group
$\Q_p^\times\!/\Q_p^{\times 2}$, which has the same $\F_2$-dimension.
Among these isomorphisms there is one which is very special.

\numeroter {\it There is a unique isomorphism\/
  $\rho_M:\Q_p^\times\!/\Q_p^{\times 2}\to G$ such that for every
  quadratic extension\/ $E$ of\/ $\Q_p$, the kernel of the composite
  map\/ $\rho_E:\Q_p^\times\to\Gal(E|\Q_p)$ is\/
  $N_{E|\Q_p}(E^\times)$.}

{\it Proof}.  We have the perfect pairing $\langle\
,\ \rangle_p:G\times(\Q_p^\times\!/\Q_p^{\times2})\to\Z^\times$ given
by $\langle\sigma,b\rangle=\sigma(\sqrt b)/\sqrt b$ for every
$\sigma\in G$ and every $\bar b\in\Q_p^\times\!/\Q_p^{\times2}$, so we
have the canonical isomorphism
$G\to\Hom(\Q_p^\times\!/\Q_p^{\times2},\Z^\times)$.  But we have just
established (??) the isomorphism
$\Q_p^\times\!/\Q_p^{\times2}\to\Hom(\Q_p^\times\!/\Q_p^{\times2},\Z^\times)$,
and hence we get an isomorphism $\rho_M:\Q_p^\times\!/\Q_p^{\times
2}\to G$, and it remains to show that $\rho_M$ has the stated
property.

This follows from the fact that when we identify these two groups
using $\rho_M$, the kummerian pairing $\langle\ ,\ \rangle_p$ gets
converted into the hilbertian pairing $(\ ,\ )_p$, and we have shown
(??) that $(a,b)_p=1$ if and only if $a$ is a norm from the extension
$\Q_p(\sqrt b)$.  \cqfd

\numeroter Note that the reciprocity isomorphism $\rho_M$ gives back
the bijection $E\mapsto\chi$ of (??).  Indeed, quadratic extensions
of $\Q_p$ correspond to quadratic characters of $G$ and hence
(applying $\rho_M^{-1}$) to quadratic characters of $\Q_p^\times$.

\numeroter Note also that $\rho_M$ respects the natural filtrations on
the groups $G$ and $\Q_p^\times\!/\Q_p^{\times2}$.  The filtration on
the latter group comes from the filtration $\cdots U_2\subset
U_1\subset\Z_p^\times\subset\Q_p^\times$.  Concretely, it is the
filtration 
$$
\langle\bar1\rangle\subset\langle\bar
u\rangle\subset\Q_p^\times\!/\Q_p^{\times2},
\quad\hbox{resp. }
\langle\bar1\rangle\subset
\langle\bar5\rangle\subset
\langle\bar5,-\bar1\rangle\subset
\Q_2^\times\!/\Q_2^{\times2}
$$ 
for $p\neq2$ and $p=2$ respectively.  For $p\neq2$, the filtration
on $G$ comes from the tower $M\mid\Q_p(\sqrt u)\mid\Q_p$ of quadratic
extensions, whereas for $p=2$ it comes from the tower
$M\mid\Q_2(\sqrt5,\sqrt{-1})\mid\Q_2(\sqrt5)\mid\Q_2$.  We shall see
later that this is the {\it ramification filtration} on $G$.

\numeroter There is a unique isomorphism
$\rho_\C:\R^\times\!/\R^{\times 2}\to\Gal(\C|\R)$ of groups, and it
has the property that $\Ker(\rho_\C)=N_{\C|\R}(\C^\times)$ analogous
to the defining property of $\rho_M$ of (??).

\numeroter What we have achieved might not seem much, but it is rare
to be able to compute explicitly, for a given galoisian extension $M$
of a field $K$, a set of elements $S\subset M$ such that $M=K(S)$, the
group $G=\Gal(M|K)$, and $\sigma(s)$ for every $\sigma\in G$ and $s\in
S$.  This is what we have done for $K=\Q_v$ and $M$ the maximal
abelian extension of exponent~$2$.  For example, when $v$ is an {\it
odd\/} prime~$p$, we may take $S=\{\sqrt u,\sqrt{-p}\}$ with
$u\in\Z_p^\times$ not a square (or equivalently $\lambda_p(u)=-1$), we
have $G=\{\sigma_1,\sigma_u,\sigma_p,\sigma_{up}\}$, and the action is
given by $\sigma_a(\sqrt b)=(a,b)_p.\sqrt b$, where $(a,b)_p$ has been
computed in (??).  A similar statement holds for $v=2$ or $v=\infty$.

\numeroter It is a minor miracle --- in my view --- that whereas
$\rho_M$ is uniquely determined by imposing the condition
$\Ker(\rho_E)=N_{E|\Q_p}(E^\times)$ on any two (resp.~three for $p=2$)
of the three (resp.~seven) quadratic extensions $E$ whose compositum
is $M$, this condition is automatically satisfied by the remaining
quadratic extension(s). In other words, $\rho_M$ is independent of the
choice of bases.

\numeroter One of the main results of the theory of {\it abelian
extensions of local fields\/} says that for every local field $K$ (of
which $\Q_p$ is the first example) and for every $n>0$, there is a
unique isomorphism $\rho:K^\times\!/K^{\times n}\to\Gal(M|K)$, where
$M$ is the maximal abelian extension of $K$ of exponent dividing $n$,
such that for every abelian extension $E\mid K$ of exponent dividing
$n$, the kernel of the resulting composite map $K^\times\to\Gal(E|K)$
is $N_{E|K}(E^\times)$ (and such that uniformisers correspond to the
canonical generator $\sigma$ of the residual extension, as opposed to
its inverse $\sigma^{-1}$).  We have proved the case $K=\Q_p$, $n=2$
(??)  (and didn't need to worry about uniformisers because an
automorphism of order~$2$ is its own inverse).

\bigbreak 
\leftline{\it Exercises}
\medskip
 
\numeroter
Define $B_k\in\Q$ in terms of the exponential series
$e^T=\exp(T)$ by the identity
$$
{T\over e^T-1}=B_0{T^0\over0!}+B_1{T^1\over1!}+\sum_{k>1}B_k{T^k\over k!},
$$
so that $B_0=1$ and $B_1=-1/2$.  Show that $B_k=0$ for every odd
$k>1$.  ({\it Hint~:} $B_0+B_1T-T/(e^T-1)$ is invariant under
$T\mapsto-T$.)  The purpose of the next few exercises is to show that
for every even integer $k>0$,  the number
$$
W_k=B_k+\sum_{l-1\,|\,k}{1\over l}
$$
(where $l$ runs through the primes such that $k\equiv0\pmod{l-1$}) is in
$\Z$ (von Staudt--Clausen, 1840).  The idea, due to Witt, is to show
that $W_k$ is in $\Z_{(p)}=\Z_p\cap\Q$ for every prime $p$.

\numeroter 
For every integer $n>0$, let $S_k(n)=1^k+2^k+\cdots+(n-1)^k$.  Show
  that $\displaystyle S_k(n) =\sum_{m\in[0,k]}{k\choose m}{B_m\over
  k+1-m}n^{k+1-m}$.  (Compare the coefficients in
$$
1+e^T+e^{2T}+\cdots+e^{(n-1)T}={e^{nT}-1\over T}{T\over e^T-1}.)
$$
Conclude that $\lim_{r\to+\infty}S_k(p^r)/p^r=B_k$ in $\Q_p$.

\numeroter  Using the fact that every
$j\in[0,p^{s+1}[$ can be uniquely written as $j=up^s+v$, where $u\in[0,p[$ and
$v\in[0,p^s[$, deduce that 
$$
{S_k(p^{s+1})\over p^{s+1}}-{S_k(p^s)\over p^s}\in\Z.
$$
(Notice that
$$
\eqalign{
S_k(p^{s+1})
&=\sum_{j\in[0,p^{s+1}[}j^k=\sum_{u\in[0,p[}\sum_{v\in[0,p^s[}(up^s+v)^k\cr
&\equiv p\big(\sum_v v^k\big)
 +kp^s\big(\sum_u u\sum_v v^{k-1}\big)\pmod{p^{2s}}\cr 
}$$
and $2\sum_uu=p(p-1)\equiv0\pmod p$, so
$
S_k(p^{s+1})\equiv pS_k(p^s)\pmod{p^{s+1}},
$
where, for $p=2$, the fact that $k$ is even has been used.)  Conclude that
$$
{S_k(p^r)\over p^r}-{S_k(p^s)\over p^s}\in\Z
$$
for any two integers $r>0$, $s>0$, and that $B_k-S_k(p)/p\in\Z_{(p)}$
(Fix $s=1$, let $r\to+\infty$, and use (??)).

\numeroter Show that
$\displaystyle
S_k(p)
\equiv\cases
{-1\pmod p&if $k\equiv0\pmod{p-1}$\cr 
\phantom{-}0\pmod p&if $k\not\equiv0\pmod{p-1}$\cr} 
$ (To see that $\sum_{j\in[1,p[}j^k\equiv0\pmod p$ when
$k\not\equiv0\pmod{p-1}$, note that, $g$ being a generator of
$\F_p^\times$, we have $g^k-1\not\equiv0\pmod p$, whereas
$$
(g^k-1)\big(\sum_{j\in[1,p[}j^k\big)
\equiv (g^k-1)\big(\sum_{t\in[0,p-1[}g^{tk}\big)
\equiv g^{(p-1)k}-1\equiv0\pmod p.)
$$
Conclude that $B_k+p^{-1}\in\Z_{(p)}$ if $k\equiv0\pmod{p-1}$ and that
$B_k\in\Z_{(p)}$ otherwise.

\numeroter In either case, the number $W_k$ (??) is in $\Z_{(p)}$ for
every $p$ (One has 
$$
W_k=\cases{
(B_k+p^{-1})+\sum_{l\neq p}l^{-1}&if $k\equiv0\pmod{p-1}$\cr 
(B_k)+\sum_l l^{-1}&if $k\not\equiv0\pmod{p-1}$\cr 
}$$ 
where $l$ runs through the primes for which $k\equiv0\pmod{l-1}$.)  Conclude
that $W_k\in\Z$ (for every even integer $k>0$).

\numeroter  Let $v$ be a place of $\Q$.  A {\it character\/} of
$\Q_v^\times$ is a continous morphism $\chi:\Q_v^\times\to\C^\times$
of groups~; we say that $\chi$ is {\it unitary\/} if
$\chi(\Q_v^\times)\subset\U$.  For every character $\chi$, the
character $\chi_1:x\mapsto\chi(x)/|\chi(x)|_\infty$ is unitary.  The
purpose of the next few exercises is to define the {\it local
constant\/} $W_v(\chi)$.

\numeroter Every character of $\R^\times$ is uniquely of the form $x\mapsto
|x|_\infty^s\lambda_\infty(x)^r$ for some $s\in\C$ and some
$r\in\F_2$ (where $\lambda_\infty$ (??) is the unique quadratic
character of $\R^\times$)~; it is unitary if and only if $s\in i\R$.
Choose a primitive $4$-th root of~$1$ in $\C$, and define the {\it
local constant\/} $W_\infty(\chi)=i^{-r}$ (so that
$W_\infty(\chi)=W_\infty(\chi_1)$).

\numeroter  Now let $p$ be a prime number and let $\chi$ be a character of
$\Q_p^\times$.  If $\chi$ is unramified (???), we put $a(\chi)=0$.
Otherwise let $n>0$ be the smallest integer such that $\chi(U_n)=1$
(???), and put $a(\chi)=n$.  The integer $a(\chi)$ is called {\it the
exponent of the conductor\/} of $\chi$ (and the ideal ${\goth
f}(\chi)=p^{a(\chi)}\Z_p$ is called {\it the conductor\/} of $\chi$).
Compute $a(\chi)$ for each of the seven quadratic characters $\chi$
(???) of $\Q_2^\times$ and each of the the three quadratic characters
$\chi$ (???) of $\Q_p^\times$ ($p\neq2$).

\numeroter Let $\chi$ be a character of $\Q_p^\times$ with associated
unitary character $\chi_1$, and choose $\gamma\in\Q_p^\times$ such
that $v_p(\gamma)=a(\chi)$, where $a(\chi)$ (???) is the exponent of
the conductor ${\goth f}(\chi)$ of $\chi$ (so that ${\goth
f}(\chi)=p^{a(\chi)}\Z_p$).  Using the notation $\langle\ \rangle_p$
from (???), define {\it the local constant\/} as
$$
W_p(\chi)
={\chi_1(\gamma)\over\sqrt p^{a(\chi)}}
 \sum_{x\in(\Z_p/{\goth f}(\chi))^\times} 
 \chi^{-1}(x)e^{2i\pi\langle x/\gamma\rangle_p}.
$$
Show that $W_p(\chi)$ does not depend on the choice of $\gamma$ and
compute it for each of the seven quadratic characters $\chi$ (???) of
$\Q_2^\times$ and each of the the three quadratic characters $\chi$
(???) of $\Q_p^\times$ ($p\neq2$).

\numeroter Let $d$ be a squarefree integer.  For each place $v$ of
$\Q$, we have the extension $\Q_v(\sqrt d)$ of $\Q_v$ of degree~$1$
or~$2$, and hence (???) a character $\chi_v$ of $\Q_v^\times$ of order
dividing~$2$.  Show that $W_v(\chi_v)=1$ for almost all $v$, and
$\prod_v W_v(\chi_v)=1$.  (Use induction on the number of prime divisors
of $d$.)

\vfill\eject

\centerline{\bf Lecture 7}
\medskip
\centerline{$|\ |$}
\bigskip

\numeroter We started with the quadratic reciprocity law (??) which
asserts that
$$
\lambda_p(-1)
=\lambda_4(p),\quad 
\lambda_p(2)
=\lambda_8(p),\quad\hbox{and\ }
\lambda_p(q)
=\lambda_q(\lambda_4(p)p)
$$ for any two distinct odd prime numbers $p$ and $q$.  We then
reformulated it as a product formula $\prod_v(a,b)_v=1$ (??), where
$v$ runs over all {\it places\/} of $\Q$, namely the prime numbers and
also the archimedean place $\infty$, and $a,b\in\Q^\times$.  Finally,
for each $v$, we've understood the factor $(a,b)_v$ (which actually
makes sense for all $a,b\in\Q_v$) in terms of a canonical isomorphism
$\rho:\Q_v^\times\!/\Q_v^{\times 2}\to\Gal(M_v|\Q_v)$, where $M_v$ is
the maximal abelian extension of exponent~$2$ of $\Q_v$ (??).  But it
is still mysterious as to why the product in (??) extends over all
prime numbers and the symbol $\infty$.  In other words, we have to
give a more intrinsic {\it definition\/} of a place of $\Q$, instead
of merely {\it declaring\/} that the primes and $\infty$ are the places of
$\Q$.

\numeroter Let $k$ be a field and denote its multiplicative group by
$k^\times$.  The multiplicative group of strictly positive reals is
denoted $\R^\times_+$~; it is a totally ordered group.  A {\it norm\/}
$|\ |$ on\/ $k$ is a homomorphism\/
$|\ |:k^\times\rightarrow\R^\times_+$ such that {\it the trinagular inequality\/}
$$
|x+y|\le |x|+|y|
$$ holds for every\/ $x,y\in k$, with the convention that\/ $|0|=0$.
A norm is called {\it trivial\/} if the image of\/ $k^\times$ is\/
$\{1\}$, {\it essential\/} otherwise.  A norm is called {\it
  unarchimedean\/} if {\it the ultrametric inequality\/}
$|x+y|\le\Sup(|x|,|y|)$ holds for all $x,y\in k$, {\it archimedean\/}
otherwise.

\numeroter If $\zeta\in k^\times$ has finite order, then $|\zeta|=1$,
for $1$ is the only element of finite order in $\R^{\times\circ}$.
Consequently, every norm on a finite field is trivial.

\numeroter We have the usual norm $|x|_\infty=\Sup(x,-x)$ on the field
$\R$~; it is archimedean.  Let $p$ be a prime number.  We have the
$p$-adic norm $|x|_p=p^{-v_p(x)}$ (??) on the field $\Q_p$~; it is
unarchimedean.

\numeroter Two norms $|\ |_1$, $|\ |_2$ on $k$ are called {\it
  equivalent\/} if there exists a real $\gamma>0$ such that
$|\ |_1=|\ |_2^\gamma$.

\numeroter {\it For every norm\/ $|\ |$ and every\/ $x,y\in k$, we have}
$$
\left|\,|x|-|y|\,\right|_\infty\le|x-y|.
$$

{\it Proof}.  As $y=x+(y-x)$, we have $|y|\le|x|+|y-x|$.  Similarly,
$|x|\le|y|+|x-y|$.  But $|y-x|=|x-y|$, hence the result. \cqfd

\numeroter {\it Let\/ $|\ |$ be an unarchimedean norm, and\/ $x,y\in
  k$ be such that\/
$|x|<|y|$.  Then\/ $|x+y|=|y|$.}

{\it Proof}.  We have $|x+y|\le\Sup(|x|,|y|)=|y|$.  On the other hand,
$y=(x+y)+(-x)$, so $|y|\le\Sup(|x+y|,|x|)$, and this $\Sup$ cannot be
$|x|$ by hypothesis. \cqfd

\numeroter {\it A norm\/ $|\ |$ is unarchimedean if and only if\/ $|\iota(n)|\le1$
for every\/ $n\in\Z$, where\/ $\iota(n)\in k$ is the image of\/ $n$.}

{\it Proof}.  If $|\ |$ is unarchimedean, then $|\iota(n)|\le1$ by
induction.  Conversely, suppose that $|\iota(n)|\le1$ for every
$n\in\Z$.  Let $x,y\in k$, and let $s>0$ be an integer.  We have
$$
\eqalign{
|x+y|^s&=\left|\sum_{r\in[0,s]}{s\choose r}x^ry^{s-r}\right|\cr
&\le\sum_{r\in[0,s]}
 \left|{s\choose r}\right|\left|x\right|^r\left|y\right|^{s-r}\cr
&\le\sum_{r\in[0,s]} \left|x\right|^r\left|y\right|^{s-r}\cr
&\le(s+1)\sup(|x|,|y|)^s.
}$$
Taking the $s$-th root and letting $s\to+\infty$, we get
$|x+y|\le\sup(|x|,|y|)$, as required. \cqfd

\numeroter {\it If the restriction of\/ $|\ |$ to a subfield is
  unarchimedean, then so is\/ $|\ |$.  Every norm on a field of
  characteristic\/ $\neq0$ is\/ unarchimedean.}

{\it Proof}.  The first statement follows from (??).  The second
statement also follows because that the restriction of\/ $|\ |$ to the
prime subfield is trivial (??). \cqfd

\numeroter A {\it valuation\/} $v$ on $k$ is a
homomorphism\ $v:k^\times\to\R$ such that
$$
v(x+y)\ge\inf(v(x),v(y))
$$ for every\/ $x,y\in k$, with the convention that\/ $v(0)=+\infty$.
A valuation\/ $v$ is said to be {\it trivial\/} if\/
$v(k^\times)=\{0\}$, of {\it height\/}~$1$ otherwise.  A height-$1$
valuation\/ $v$ is called {\it discrete\/} if the subgroup\/
$v(k^\times)\subset\R$ is discrete. A discrete valuation\/ $v$ is said
to be {\it normalised\/} if\/ $v(k^\times)=\Z$.
    
\numeroter Notice that if $v(x)<v(y)$, then $v(x+y)=v(x)$.  The proof
is the same as that for (??).

\numeroter For every prime $p$, we have the $p$-adic valuation $v_p$
(??) on the field $\Q_p$~; it is discrete and normalised.

\numeroter Two valuations $v_1$, $v_2$ are called {\it equivalent\/}
if there exists a real $\gamma>0$ such that $v_1=\gamma v_2$.

\numeroter If $v$ is a height-$1$ valuation on $k$, then
$|x|_v=\exp(-v(x))$ is an essential unarchimedean norm, and
conversely, if $|\ |$ is an essential unarchimedean norm, then
$v_{|\ |}(x)=-\log|x|$ is a height-$1$ valuation.  Equivalent
valuations correspond to equivalent (unarchimedean) norms.

\numeroter (Ostrowski, 1918) {\it Let\/ $|\phantom{x}|$ be a norm on
  the field\/ $\Q$ of rational numbers.  Then\/ $|\phantom{x}|$ is
  either trivial, or equivalent to the archimedean norm\/
  $|\phantom{x}|_\infty$, or equivalent to the $p$-adic norm\/
  $|\phantom{x}|_p\;$ for some prime\/ $p$.}

{\it Proof\/} (Artin, 1932). Clearly $|\phantom{x}|$ is trivial if
$|p|=1$ for every prime $p$ (because the group $\Q^\times$ is
generated by the set of primes and $-1$).  Assume that $|p|\neq1$ for
some prime $p$. We shall show that if $|p|>1$, then $|l|>1$ for every
prime $l$ and that $|\phantom{x}|$ is equivalent to
$|\phantom{x}|_\infty$.  On the other hand, if $|p|<1$, then $|l|=1$
for every prime $l\neq p$, and that $|\phantom{x}|$ is equivalent to
$|\phantom{x}|_p$.

For the time being, let $p$ and $l$ be any two integers $>1$, and
write $p$ as 
$$
p=a_0+a_1l+a_2l^2+\cdots+a_nl^n%
$$ 
in base~$l$, with digits $a_i\in[0,l[$ and $l^n\le p$,
    i.e.,~$n\le\alpha$, with $\displaystyle \alpha={\log p\over\log
      l}$.  The triangular inequality gives
$$\eqalign{
|p|&\le |a_0|+|a_1||l|+|a_2||l|^2+\cdots+|a_n||l|^n\cr
&\le(|a_0|+|a_1|+|a_2|+\cdots+|a_n|)\,\sup(1,|l|^n)\cr
&\le (1+n)d\,\sup(1,|l|^n)
\qquad(\hbox{with\ } d=\sup(|0|, |1|, \ldots, |l-1|)\cr
&\le (1+{\alpha})d\,\sup(1,|l|^{\alpha})
\qquad\left(\hbox{since\ } n\le{\alpha}\right).\cr
}
$$
Replace $p$ by $p^s$ and extract the $s$-th root to get
$$
|p|\le \left(1+s{\alpha}\right)^{1/s}d^{1/ s}
\sup\left(1,|l|^{\alpha}\right),
$$
so that we obtain the estimate
$
|p|\le\sup\left(1,|l|^{\alpha}\right)
$
upon letting $s\to+\infty$.

Now suppose that $p$ and $l$ are prime numbers and that $|p|>1$.  As
$1<|p|\le\sup\left(1,|l|^{\alpha}\right)$, we see that $|l|^{\alpha}>1$ and
hence $|l|>1$ and $|p|\le|l|^{\alpha}$. Interchanging the role of $p$ and $l$,
we deduce $|p|=|l|^{\alpha}$.

Defining $\gamma$ ($>0$) by the equation $|p|=|p|_\infty^\gamma$ for the fixed
prime $p$, we see that $|l|=|l|_\infty^\gamma$ for every prime $l$,
and hence $|x|=|x|_\infty^\gamma$ for every $x\in\Q^\times$,
i.e.,~$|\phantom{x}|$ is equivalent to the archimedean absolute value
$|\phantom{x}|_\infty$.

Finally, assume that $|p|<1$ for some prime $p$.  
We have already seen that then $|l|\le1$ for every prime $l$, and
hence $|n|\le1$ for every $n\in\Z$.  Let us show that $|l|=1$ for
every prime $l\neq p$.  For every integer $s>0$, writing
$1=a_sp+b_sl^s$ ($a_s,b_s\in\Z$), we have
$$
1=|1|\le|a_s|\,|p|+|b_s|\,|l|^s\le|p|+|l|^s,
$$
i.e.,~$|l|^s\ge1-|p|$.  This is possible for all integers $s>0$ only if
$|l|=1$.  This shows that $|\phantom{x}|$ is equivalent to
$|\phantom{x}|_p$. \cqfd

\numeroter A {\it place\/} of $k$ is an equivalent class of essential
norms on $k$.  The foregoing theorem asserts that places of $\Q$
correspond naturally to the set $\bar P$ of prime numbers (the
unarchimedean places) together with $\infty$, the archimedean place.
The following theorem, called the product formula, is another piece of
evidence for $\bar P$ being the set of all places of $\Q$.

\numeroter {\it For every $x\in\Q^\times$, one has\/ $|x|_v=1$ for
  almost all\/ $v\in\bar P$ and\/ $\prod_{v\in\bar P}|x|_v=1$.}

{\it Proof}.  Indeed, it is sufficient to verify this for $x=-1$ and
for $x=p$ (where $p$ is a prime number).  Note that for the product
formula to hold, we have to normalise the norms on $\Q$
suitably. \cqfd

\numeroter {\it If the image of a norm\/
  $|\ |:k^\times\to\R^\times_+$ is discrete, then\/  $|\ |$ is
  unarchimedean.}

{\it Proof}.  If $k$ has characteristic~$\neq0$, then every norm is
unarchimedean (??), so there is nothing to prove.  If $k$ has
characteristic~$0$, then the restriction of $|\ |$ to $\Q$ has
discrete image by hypothesis, so must be trivial or equivalent to
$|\ |_p$ for some prime $p$ (??).  Hence $|\ |$ is unarchimedean
(??). \cqfd

\medbreak
\leftline{\it  Exercises}  
\medskip

\numeroter Let $k$ be a field and put $K=k(T)$.  Recall that the group
$K^\times\!/k^\times$ is the free commutative group on the set $P_K$
of monic irreducible polynomials $f$ in $k[T]$.  For each $f\in P_K$,
let $v_f:K^\times\rightarrow\Z$ be the unique homomorphism which is
trivial on $k^\times$, sends $f$ to $1$, and sends every other element
of $P_K$ to $0$.  Check that $v_f$ is a discrete valuation.  Also, the
map $v_\infty:K^\times\rightarrow\Z$ which sends $a$ to $-\deg(a)$ is
a discrete valuation, trivial on $k^\times$.  The discrete valuations
$v_\infty$, $v_f$ (for varying $f\in P_K$) are mutually inequivalent.

\numeroter Up to equivalence, the only height-$1$ valuations on\/
$k(T)$, trivial on\/ $k$, are\/ $v_\infty$ and the\/ $v_f$, one for
each\/ $f\in P_K$.  (Let $v$ be a height-$1$ valuation on
$k(T)$, trivial on $k$.  We will show that if $v(T)<0$, then $v$ is
equivalent to $v_\infty$, whereas if $v(T)\ge0$, then there is a
unique $f\in P_K$ with $v(f)>0$ and $v$ is equivalent to $v_f$.

Suppose that $v(T)<0$.  It is sufficient to show that
$v(f)=v(T)\deg(f)$ for every $f\in P_K$.  This is clearly true for
$f=T$.  For any other $f$, write
$f=T^{n}(1+\alpha_1T^{-1}+\cdots+\alpha_nT^{-n})$ (with $n=\deg(f)$ and
$\alpha_i\in k$, at least one of them $\neq0$).  As $v(1)=0$ and
$v(\alpha_1T^{-1}+\cdots+\alpha_nT^{-n})>0$, we have
$v(1+\alpha_1T^{-1}+\cdots+\alpha_nT^{-n})=0$, and, finally,
$v(f)=v(T)\deg(f)$, i.e.~$v$ is equivalent to $v_\infty$.

Suppose now that $v(T)\ge0$~; then $v(a)\ge0$ for all $a\in k[T]$.  If
we had $v(f)=0$ for every $f\in P_K$, the valuation $v$ would be
trivial, not of height~$1$.  Pick $p\in P_K$ for which $v(p)>0$.
For every $q\neq p$ in $P_K$, write $1=ap+bq$ ($a,b\in k[T]$).  We
have
$$
0=v(1)\ge\inf\left(v(a)+v(p),v(b)+v(q)\right)
\ge\inf\left(v(p),v(q)\right)\ge0,
$$
which is possible only if $v(q)=0$.  It follows that $v$ is equivalent
to $v_p$.  It is instructive to compare this proof with Artin's proof classifying
absolute values on $\Q$.)

\numeroter Up to equivalence, the only norms on\/ $k(T)$, trivial on\/
$k$, are\/ $|\phantom{x}|_\infty$ and the\/ $|\phantom{x}|_f$, one for
each\/ $f\in P_K$. (Put $\bar P_K=P_K\cup\{\infty\}$ and
$\deg(\infty)=1$.  For each $p\in\bar P_K$, define the norm
$|x|_p=e^{-v_p(x)\deg(p)}$.)

\numeroter 
For every\/ $x\in k(T)^\times$, one has\/ $|x|_p=1$ for almost all\/
$p\in\bar P_K$ and\/ $\prod_{p\in\bar P_K}|x|_p=1$. ( This is clearly
true for $x\in k^\times$, so it is sufficient to check this for $x\in
P_K$.  We have $v_\infty(x)=-\deg(x)$, $v_x(x)=1$, and $v_q(x)=0$ for
every $q\neq x$ in $P_K$, which gives the ``sum formula''
$\sum_{p\in\bar P_K}\deg(p)v_p(x)=0$.  The result follows from this
upon exponentiating.  Note that here too, as in the case of $\Q$
earlier, it is necessary to normalise the absolute values suitably for
the product formula to hold.)

\vfill\eject

\centerline{\bf Lecture 8}
\medskip
\centerline{$ax^2+by^2=1$}
\bigskip

\numeroter Information gleaned locally at each place $v$ of $\Q$
(namely in the fields $\Q_v$) can often be pieced together to get a
global result (about the field $\Q$).  The simplest such example is
the main result (???) of this lecture.  It is further evidence for the
fact that we have found the right notion of a place of $\Q$.

\numeroter We first need a simple lemma for the proof.  Let $k$ be any
field, and $a,b,c\in k^\times$, $d\in k$ such that $d^2-a=bc$. Let $S$
be the set of $(x,y,s)\neq(0,0,0)$ in $k^3$ such that $ax^2+by^2=s^2$
and let $T$ be the set of $(w,z,t)\neq(0,0,0)$ in $k^3$ such that
$aw^2+cz^2=t^2$.  It is easy to see that if $(x,y,s)\in S$, then
$(dx+s,by,ax+ds)\in T$ and if $(w,z,t)\in T$, then
$(dw-t,cz,-aw+dt)\in S$.

\numeroter {\it The maps\/ $S\to T$ and\/ $T\to S$ defined above are bijections,
reciprocal to each other.}  \cqfd

\numeroter {\it Let\/ $a,b\in\Q^\times$.  For there to exist
  $x,y\in\Q$ such that\/ $ax^2+by^2=1$, it is necessary and sufficient
  that\/ $(a,b)_v=1$ for every place\/ $v$ of\/ $\Q$.}

\numeroter The local conditions ($(a,b)_v=1$ for
every~$v$) are clearly necessary (??).  Let us show that they are also
sufficient.  So assume that $(a,b)_v=1$ for every place $v$.

\numeroter The existence of $x,y\in\Q$ such that\/ $ax^2+by^2=1$ is
unaffected if we replace $a,b$ by $ac^2,bd^2$ for some
$c,d\in\Q^\times$~; also, $(ac^2,bd^2)_v=(a,b)_v$ for every~$v$.  We
may therefore assume that $a$ and $b$ are squarefree integers, and
proceed by induction on $|a|_\infty+|b|_\infty$.

\numeroter If either $a=1$ or $b=1$, then we have the solution
$(x,y)=(1,0)$ or $(x,y)=(0,1)$, as the case may be.  If $|a|_\infty+|b|_\infty=2$,
then we must have $|a|_\infty=1$ and $|b|_\infty=1$, and either $a=1$ or $b=1$
(since $(a,b)_\infty=1$), hence there is a solution $(x,y)$, as we
have just seen.

\numeroter Suppose that $|a|_\infty+|b|_\infty>2$, and assume (up to
interchaning $a$ and $b$) that $|a|_\infty\le|b|_\infty$.  Let us first show that
the local conditions force the existence of an $d\in\Z$ such that
$a\equiv d^2\pmod b$.  As $b$ is squarefree, it is sufficient to show
that $\lambda_p(a)=1$ (unless $a\equiv0\pmod p$) for every {\it odd\/}
prime divisor $p|b$.  But $\lambda_p(a)=(a,b)_p=1$, where the first
equality holds because $v_p(b)=1$ (??).  So the existence of $d$ is
guaranteed, and we may further assume that $d\in[0,|b|_\infty/2]$.

\numeroter Put $d^2-a=bc$ for some $c\in\Z$.  If $c=0$, then we have
the solution $(x,y)=(d^{-1},0)$.  Suppose that $c\neq0$~; then
$$
|c|_\infty=\left|d^2-a\over b\right|_\infty\le
\left|d^2\over b\right|_\infty+\left|a\over b\right|_\infty\le
{|b|_\infty\over4}+1<|b|_\infty,
$$
where the last inequality holds since $|b|_\infty>1$.  

\numeroter But we have seen (??) that, over any field $k$, the
existence of a solution $(x,y)\in k^2$ of $ax^2+by^2=1$ is equivalent to
the existence of a solution $(w,z)\in k^2$ of $aw^2+cz^2=1$, for any given
$a,b,c\in k^\times$ and $d\in k$ such that $d^2-a=bc$.

\numeroter Let's return to our $a,b,c,d$ from (??), and write $c=ef^2$
for some $e,f\in\Z$ of which $e$ is squarefree, and note that
$|e|_\infty<|b|_\infty$.  The proof is therefore over by the inductive
hypothesis, since $|a|_\infty+|e|_\infty<|a|_\infty+|b|_\infty$.
Indeed, as there are local solutions for $ax^2+by^2=1$ at every place
of $\Q$ by hypothesis, there are local solutions everywhere for
$aw^2+cz^2=1$ (applying (??) to $k=\Q_v$) and hence for $aw^2+et^2=1$,
and therefore a global solution for $aw^2+et^2=1$ (the inductive step,
to apply which we changed $c$ into the squarefree~$e$) and hence for
$aw^2+cz^2=1$, and therefore a global solution for $ax^2+by^2=1$
(applying (??)  with $k=\Q$).  \cqfd

\numeroter  {\it Let\/ $a,b,c\in\Q_v^\times$ and\/ $r\in\Q_v$ be such that\/ $r^2-a=bc$.
Then\/ $(a,b)_v=(a,c)_v$.  In particualr, $(r^2-a,a)_v=1$.}

{\it Proof}.  ??? Notice that for $r=0$ and $r=1$, we recover the
relations $(-a,a)_v=1$ and $(1-a,a)_v=1$ proved earlier.  \cqfd

\numeroter {\it Let\/ $a,b\in\Q$.  Then $a\in N_{\Q(\sqrt
    b)|\Q}(\Q(\sqrt b)^\times)$ if and only if\/ $a\in N_{\Q_v(\sqrt
    b)|\Q_v}(\Q_v(\sqrt b)^\times)$ for every place\/ $v$ of\/ $\Q$.}

{\it Proof}.  ??? This is expressed by saying that $a$ is a norm from the
  extension $\Q(\sqrt b)$ if and only if it is everywhere locally a
  norm. \cqfd

\bigbreak 
\leftline{\it Exercises}
\medskip

\numeroter (Legendre) Let $a,b,c\in\Z$ be integers such that $abc\neq0$ is
squarefree, and suppose that $a,b,c$ are not all three of the same
sign.  Then $ax^2+by^2+cz^2=0$ has a solution $(x,y,z)\neq(0,0,0)$ in
$\Z^3$ if and only if $-bc$, $-ca$, $-ab$ are squares $\pmod a$,
$\pmod b$, $\pmod c$ respectively.

\numeroter The law of quadratic reciprocity (??) was not used in the
proof of Legendre's theorem (??).  Criticise the following purported
proof of the said law from this theorem~: We consider eight cases
according to the signs of $\lambda_4(p)$, $\lambda_4(q)$ and
$\lambda_q(p)$ (as Gau\ss\ did in his first proof).  In each case, we
apply (??) to a suitable triple $(a,b,c)$ such that $a\equiv b\equiv
c\equiv1\pmod4$ (conditions which force
$ax^2+by^2+cz^2\not\equiv0\pmod4$, cf.~(??))  to determine
$\lambda_p(q)$.  For example, here are some cases : \par\noindent

{\it i\/}) If $\lambda_4(p)=1$, $\lambda_4(q)=-1$ and $\lambda_q(p)=-1$, we
consider the triple $(a,b,c)=(1,p,-q)$ and conclude that
$\lambda_p(q)=-1$, as required.

{\it ii\/}) Similarly, if $\lambda_4(p)=-1$, $\lambda_4(q)=-1$ and $\lambda_q(p)=1$, we
take $(a,b,c)=(1,-p,-q)$.

{\it iii\/}) Now consider the case $\lambda_4(p)=-1$, $\lambda_4(q)=-1$,
$\lambda_q(p)=-1$, let $l$ be an auxillary prime such that
$\lambda_4(l)=1, \lambda_p(l)=-1$, $\lambda_q(l)=-1$, and take
$(a,b,c)=(l,-p,-q)$ to conclude that $\lambda_p(q)=1$.

{\it iv\/}) In case $\lambda_4(p)=1$ and $\lambda_4(q)=1$, consider an
auxillary prime $l$ such that $\lambda_4(l)=-1$, $\lambda_q(l)=1$,
$\lambda_l(p)=-1$, and take $(a,b,c)=(p,q,-l)$.

\vfill\eject

\bigbreak
\unvbox\bibbox

\bye